\documentclass[12pt]{article}
\usepackage{amssymb}
\usepackage{amsmath}
\usepackage{array}
\usepackage[T2A]{fontenc}
\usepackage[cp1251]{inputenc}
\usepackage[pdftex]{graphics}
\renewcommand{\S}{\mathhexbox278}

\renewcommand{\S}{\mathhexbox278}

\DeclareMathOperator{\RRe}{Re} \DeclareMathOperator{\IIm}{Im}

\DeclareMathOperator{\vth}{\vartheta}
\DeclareMathOperator{\vep}{\varepsilon}
\DeclareMathOperator{\vk}{\varkappa}

\newcommand{\D}{\displaystyle}

\newcommand{\prsum}{\mathop{{\sum}'}}

\renewcommand{\le}{\operatorname{\leqslant}}
\renewcommand{\ge}{\operatorname{\geqslant}}

\multlinegap=0em \DeclareFontFamily{T1}{msb}{}
\DeclareFontShape{T1}{msb}{m}{ol}{<5> <6> <7> <8> <9> gen * msbm
<10> <10.95> <12> <14.4> <17.28> <20.74> <24.88> msbm10}{}
\DeclareSymbolFont{AMSb}{T1}{msb}{m}{ol} \multlinegap=0em
\begin{document}

\begin{center}
{\rmfamily\bfseries\normalsize On the Gram's Law in the Theory of
Riemann Zeta Function\footnote{This research was supported by the
Programme of the President of the Russian Federation `Young
Candidates of the Russian Federation' (grant no. МК-4052.2009.1).}}
\end{center}

\begin{center}
{\rmfamily\bfseries\normalsize M.A. Korolev}
\end{center}

\vspace{0.5cm}

\fontsize{11}{12pt}\selectfont

\textbf{Abstract}. Some statements concerning the distribution of
imaginary parts of zeros of the Riemann zeta\,-function are
established. These assertions are connected with so\,-called `Gram's
law' or `Gram's rule'. In particular, we give a proof of several
Selberg's formulae stated him without proof in his paper `The Zeta
Function and the Riemann Hypothesis' (1946), and some of their
equivalents.

\fontsize{12}{15pt}\selectfont

\begin{center}
\textbf{Introduction}
\end{center}

In the present paper, the author continues his studies begun in
\cite{Korolev_2010} and connected with the distribution of ordinates
of zeros of the Riemann zeta\,-function $\zeta(s)$ and with
so\,-called Gram's law.

Let's remind some basic notions and definitions.

\vspace{0.2cm}

\textbf{Definition 1.} For positive $t$, we denote by $\vth(t)$ the
increment of the argument of the function
$\pi^{-s/2}\Gamma\bigl(\frac{s}{2}\bigr)$ along the segment with the
end\,-points $s = \tfrac{1}{2}$ and $s = \tfrac{1}{2}+it$.

\vspace{0.2cm}

On can prove (see, for example, \cite{Siegel_1943}) that the
asymptotic expansion
\[
\vth(t)\,\sim\,\frac{t}{2}\ln\frac{t}{2\pi}\;-
\;\frac{t}{2}\;-\;\frac{\pi}{8}\;+\;\sum\limits_{n =
1}^{+\infty}\frac{2^{2n-1}-1}{2^{2n}\;2n(2n-1)}\;(-1)^{n+1}B_{2n}\;t^{-(2n-1)},
\]
holds, as $t$ grows. Here the $B_{2n}$ are Bernoulli numbers: $B_{2}
= \tfrac{1}{6}$, $B_{4} = -\tfrac{1}{30}$, $B_{6} = \tfrac{1}{42}$,
$B_{8} = -\tfrac{1}{30}$, $B_{10} = \tfrac{5}{66}$ and so on.

The function $\vth(t)$ is presented in the expression for $N(t)$ --
the number of zeros of $\zeta(s)$ in the strip $0<\IIm s\le t$. This
expression is called by Riemann - von Mangoldt formula and has the
form
\[
N(t)\,=\,\frac{1}{\pi}\,\vth(t)\,+\,1\,+\,S(t).
\]
Here $S(t) = \pi^{-1}\arg{\zeta\bigl(\tfrac{1}{2}+it\bigr)}$ denotes
the argument of the Riemann zeta function on the critical line. For
basic properties of $S(t)$, see the survey
\cite{Karatsuba_Korolev_2005}.

\vspace{0.2cm}

\textbf{Definition 2.} For any $n\ge 0$, the \emph{Gram point}
$t_{n}$ is defined as the unique root of the equation
$\vth(t_{n})\,=\,\pi\!\cdot\!(n-1)$.

\vspace{0.2cm}

Suppose $0<\gamma_{1}<\gamma_{2}<\ldots \le \gamma_{n}\le
\gamma_{n+1}\le \ldots$ are positive ordinates of zeros of
$\zeta(s)$ numbering in ascending order (if several zeros have the
same ordinates, we numerate them in arbitrary way). First
researchers of zeros of the Riemann zeta function observed that for
all the values $n$ not very large, with the exception of small part
of cases, the ordinate $\gamma_{n}$ belongs to the interval $G_{n} =
(t_{n-1}, t_{n}]$. Later, this phenomenon was called as `Gram's law'
or `Gram's rule'.

If someone wants to formulate the general rule based on the small
number of examples, he has some freedom. This is the reason why the
words `Gram's rule' (and `Gram's law') in some papers have the sense
different from the above. Thus, one can say that the interval
$G_{n}$ satisfies the Gram's rule if $G_{n}$ contains exactly one
zero of function $\zeta\bigl(\tfrac{1}{2}+it\bigr)$, and the number
of this zero doesn't matter. This version of Gram's rule is
contained in the papers of J.~I.~Hutchinson \cite{Hutchinson_1925}
and E.~C.~Titchmarsh \cite{Titchmarsh_1935}, and from this point of
view Gram's rule is studied by T.~S.~Trudgian \cite{Trudgian_2009}.

In order to characterize the degree of deviation from the Gram's
rule, Titchmarsh considered \hfill the \hfill fractions \hfill
$\frac{c_{n}-t_{n}}{t_{n+1}-t_{n}}$ \hfill where \hfill $c_{n}$
\hfill denotes \hfill the \hfill positive \hfill zeros \hfill of \hfill the \hfill function \\
$\zeta\bigl(\tfrac{1}{2}+it\bigr)$ numerated in ascending order, and
proved their unboundedness. If the interval $G_{m}$ contains the
zero $c_{n}$ then such fraction is close to the difference $m-n$.
Thus, Titchmarsh's theorem implies that the number $n$ of zero lying
on the critical line can differs from the number of corresponding
Gram's interval by arbitrary large value. In particular, this means
that there exist infinitely many exceptions from Gram's rule.

In the part IV of his report `The Zeta\,-Function and the Riemann
Hypothesis' \cite{Selberg_1946a}, Selberg formulated Titchmarsh's
theorem in a following way:
\begin{equation}
\varliminf_{n\to +\infty}\Delta_{n}\,=\,-\infty, \quad
\varlimsup_{n\to +\infty}\Delta_{n}\,=\,+\infty.\label{Lab3}
\end{equation}
For given $n\ge 1$, the quantity $\Delta_{n}$ was defined as $m-n$
if the following inequalities hold: $t_{m-1}<\gamma_{n}\le t_{m}$.
In 1940's, Selberg established that the formulae
\begin{equation}
\sum\limits_{n =
1}^{N}\Delta_{n}^{2k}\,=\,\frac{(2k)!}{k!}\,\frac{N}{(2\pi)^{k}}\,(\ln\ln
N)^{k}\,+\,O\bigl(N(\ln\ln N)^{k-1/2}\bigr),\label{Lab4}
\end{equation}
\begin{equation}
\sum\limits_{n = 1}^{N}\Delta_{n}^{2k-1}\,=\,O\bigl(N(\ln\ln
N)^{k-1}\bigr)\label{Lab5}
\end{equation}
are valid for any fixed integer $k\ge 1$ and stated them (without
proof) in \cite{Selberg_1946a}. Then he posed the following
conjecture: `\emph{It is probably true}, \emph{though I have not
been able to prove it rigorously}, \emph{that} $\sqrt{\log\log
n\mathstrut}$ \emph{is the} ``\emph{normal}'' \emph{order of
magnitude of} $\Delta_{n}$ \emph{in the sense that if} $\Phi(n)$
\emph{is a positive function of} $n$ \emph{which tends to infinity
with} $n$, \emph{the inequalities}
\[
\frac{\sqrt{\log\log
n}}{\Phi(n)}\;<\;|\Delta_{n}|\;<\;\Phi(n)\sqrt{\log\log n}
\]
\emph{hold true for almost all} $n$. \emph{In particular this should
imply that} $\gamma_{\nu}$ ``\emph{almost never}'' \emph{lies in the
interval} $(t_{\nu-1}, t_{\nu})$. \emph{It is the first inequality
which is difficult point}, \emph{and which I have not been able to
prove completely}.'

Thus the cases considered earlier as `exceptions' are normal, and
the cases called earlier as a `rule' occur very rarely.

There are some serious reasons to think that in his definition of
$\Delta_{n}$ Selberg denoted by $\gamma_{n}$ the ordinates of all
zeros of $\zeta(s)$ and not only the ordinates of zeros lying on the
critical line, and considered the quantities
\begin{equation}
q_{n}\,=\,\frac{\gamma_{n} - t_{n}}{t_{n+1}-t_{n}}\label{Lab2}
\end{equation}
instead of Titchmarsh's fractions. Indeed, if one considers only the
zeros on the critical line in the definition of $\Delta_{n}$, then
the formulae (\ref{Lab4}) and (\ref{Lab5}) imply some very sharp
statements concerning the distribution of zeros of $\zeta(s)$. These
statements are close to that the Riemann hypothesis asserts and, in
any case, they are much more powerful that all recent results about
a part of zeros of zeta\,-function lying on the critical line. This
is the reason why we follow the below definitions.

\vspace{0.2cm}

\textbf{Definition 3.} We say that the ordinate $\gamma_{n}$
satisfies to the Gram's law if the following inequalities hold:
$t_{n-1}<\gamma_{n}\le t_{n}$.

\vspace{0.2cm}

This means that in the case when the interval $G_{n}$ contains all
the ordinates $\gamma_{n-s},\!\ldots,$ $\gamma_{n-1}, \gamma_{n},
\gamma_{n+1},$ $\ldots, \gamma_{n+r}$, all of these ordinates,
except $\gamma_{n}$, do not obey the Gram's law.

\vspace{0.2cm}

\textbf{Definition 4.} Suppose for given $n\ge 1$ the following
inequalities hold:
\begin{equation}
t_{m-1}<\gamma_{n}\le t_{m}.\label{Lab1}
\end{equation}
Then we set $\Delta_{n} = m-n$.

\vspace{0.2cm}

Now it's clear that the ordinate $\gamma_{n}$ satisfies the Gram's
law if and only if $\Delta_{n} = 0$.

In \cite{Korolev_2010} the author proved the part of Selberg
conjecture that asserts that the ordinate $\gamma_{n}$ `almost
never' lies in the interval $G_{n} = (t_{n-1}, t_{n}]$. The proof is
based on the properties of the sequence $\Delta(n)$ defined as
follows.

\vspace{0.2cm}

\textbf{Definition 5.} Suppose for given Gram point $t_{n}$ the
following inequalities hold:
\begin{equation}
\gamma_{m}\le t_{n} < \gamma_{m+1}.\label{Lab6}
\end{equation}
Then we set $\Delta(n) = m-n$.

\vspace{0.2cm}

It appears that for given interval the number of the indexes $n$
such that $\Delta_{n}$ equals to zero, is very close to the number
of $\Delta(n)$ with the same property. Next, the study of the
sequence $\Delta(n)$ is much more simple than the study of
$\Delta_{n}$. The reason is that both (\ref{Lab6}) and Riemann\,-von
Mangoldt formula imply the key equality $\Delta(n) = S(t_{n})$. This
relation reduces the initial problem to some problem concerning the
distribution of values of the argument of the Riemann zeta function.
The behavior of the function $S(t)$ on the `regular' sequence of
Gram points is studied by methods belonging to A.~Selberg
\cite{Selberg_1946b} and A.~Ghosh \cite{Ghosh_1983}.

After the paper \cite{Korolev_2010} was published, the author
familiarized with the book of collected papers of A.Selberg issued
in 1989\,-1991. The text of Selberg's paper mentioned above in the
first volume (see \cite[с. 341\,-355]{Selberg_1989}) was provided by
the following remark: `\ldots \emph{it follows from these equations}
(i.e. from the relations (\ref{Lab4}) and (\ref{Lab5}) -- M.K.)
\emph{by standard theory that the quantity}
$\Delta_{n}/\sqrt{\log\log{n}\mathstrut}$ \emph{has a simple
Gaussian distribution. In particular this answers in the affirmative
question raised here} (i.e. the Selberg's conjecture -- M.K.).
\emph{In 1946 I did not know that these results about the moments
of} $\Delta_{n}$ \emph{allow one to determine this distribution
function.}' This remark needs careful explanations.

Indeed, before the book \cite{Selberg_1989} was printed, the
articles of A.~G.~Postnikov \cite{Postnikov_1966} and
V.~F.~Gaposchkin \cite{Gaposhkin_1968} were already published. In
these papers, the procedure of deriving of the distribution function
of some random quantity from the formulae for even and odd moments
of this quantity, was introduced. Next, in the later papers of
A.~Fujii \cite{Fujii_1976} and A.~Ghosh \cite{Ghosh_1983} such
arguments are applied to some problems in the theory of Riemann zeta
function and Dirichlet's $L$\,-functions (for the Ghosh method, see,
for example, the survey \cite{Karatsuba_Korolev_2006}).

But there are some uncertainties with the justification of the
formulae (\ref{Lab4}) and (\ref{Lab5}). As far as the author found
out, Selberd had never published his proof. However, the attempt of
such proof is presented in the paper of A.Fujii mentioned above,
where the similar problem for Dirichlet's $L$\,-functions is
considered. Let's consider this attempt more closely.

Suppose $q\ge 5$ is fixed and denote by $\chi_{1}$, $\chi_{2}$
different primitive characters modulo $q$. Let's enumerate the
positive ordinates of complex zeros of the functions $L(s,\chi_{1})$
and $L(s,\chi_{2})$ in ascending order:
\begin{align*}
& 0 < \gamma_{1}(\chi_{1}) \le \gamma_{2}(\chi_{1})\le \ldots \le
\gamma_{n}(\chi_{1})\le \gamma_{n+1}(\chi_{1})\le \ldots,\\
& 0 < \gamma_{1}(\chi_{2}) \le \gamma_{2}(\chi_{2})\le \ldots \le
\gamma_{n}(\chi_{2})\le \gamma_{n+1}(\chi_{2})\le \ldots.
\end{align*}

Next, suppose for given $n$ the following inequalities hold:
\[
\gamma_{m}(\chi_{1})<\gamma_{n}(\chi_{2})\le \gamma_{m+1}(\chi_{1}).
\]
Then we define $\Delta_{n}(\chi_{1},\chi_{2}) = n-m$. The question
is: how often the difference $\Delta_{n}(\chi_{1},\chi_{2})$
vanishes as the ordinate $\gamma_{n}(\chi_{2})$ varies in the
interval $(T, T+H]$ which is supposed to be long enough?

In order to solve this problem, Fujii proved that the discrete
random quantity with the values
$\Delta_{n}(\chi_{1},\chi_{2})/\sqrt{\ln\ln n\mathstrut}$ has the
distribution which tends to Gaussian distribution as $T$ grows. This
implies that the values of $\Delta_{n}(\chi_{1},\chi_{2})$ differ
from zero for `almost all' $n$. However, the proof of these facts is
based on the consideration of the integral
\[
j\,=\,\int_{T}^{T+H}\bigl(S(t,\chi_{1})\,-\,S(t,\chi_{2})\bigr)^{2k}dS(t,\chi_{2})
\]
and, in particular, on the estimation
\begin{equation}
j\,=\,O\bigl((\ln T)^{2k+1}\bigr).\label{Lab7}
\end{equation}
This estimation plays the key role in Fujii's arguments. Here
$S(t,\chi) = \pi^{-1}\!\arg{L\bigl(\tfrac{1}{2}+it,\!\chi\bigr)}$
denotes the argument of Dirichlet's $L$\,-function $L(s,\chi)$ on
the critical line. The proof of (\ref{Lab7}) is missed.

This proof seems not satisfactory because the integral $j$ as
Stiltjes integral does not exist. Indeed, the points of
discontinuity of piecewise smooth function $S(t,\chi_{2})$ under the
differential coincides with the ordinates $\gamma_{n}(\chi_{2})$.
The function $S(t,\chi_{1})-S(t,\chi_{2})$ in the integrand has
discontinuities at these points, too. In such case it is easy to
construct two sequences of the integral sums for $j$ which tend to
different limits (see, for example, \cite[$n^{\circ}$.
584]{Fihtengolz_1966})).

These arguments, applied formally to the problem of calculating of
even moments of $\Delta_{n}$, lead to the integrals
\[
\int_{T}^{T+H}S^{2k}(t)dS(t),\quad k = 1,2,3,\ldots,
\]
which do not exist, too.

Thus the problem of reconstruction of Selberg's proof of
(\ref{Lab4}) and (\ref{Lab5}) is still open. In the below, we make
such attempt. Now let's consider the structure of the paper.

This article consist of three parts. The first part (\S 1) deals
with the order of growth of $\Delta_{n}$ as $n\to +\infty$ (theorem
1). Next, an $\Omega$\,-estimations, which improve the assertions
(\ref{Lab3}), are proved here (theorem 2). In essence, all these
statements follow from the famous $O$\,-- and $\Omega$\,-- theorems
for the function $S(t)$.

Lemma 2 in the basic assertion of the second part of the paper (\S
2). This lemma deals with the number of solutions of the inequality
\[
a<\Delta_{n}\le b,
\]
for given integer $a$ and $b$. It appears that this number is very
close to the number of solutions of the inequality
\[
-(b+1)<\Delta(n)\le -(a+1)
\]
(the domains of $n$ are the same in both inequalities, of course).
This fact implies the closure of the distributions of
 $\Delta_{n}$ and $\Delta(n)$ (theorem 3) and the closure of their
moments of any degree. Thus, both the Selberg's conjecture and the
formulae (\ref{Lab4}),(\ref{Lab5}) follow from the corresponding
theorems of \cite{Korolev_2010} (see lemmas 3\,-5 of the present
paper).

Finally, in the last, third part some equivalents of Selberg's
conjecture are proved. In particular, the problem of distribution of
the differences $\gamma_{n}-t_{n}$ is considered (theorems 7\,-10).
If the ordinate $\gamma_{n}$ obeys the Gram's law then the order of
such difference does not exceed the quantity
\begin{equation}
t_{n+1}\,-\,t_{n}\,\sim\,\frac{2\pi}{\ln n}.\label{Lab8}
\end{equation}
However, the quantities $|\gamma_{n}-t_{n}|$ are much larger then
(\ref{Lab8}) in the mean and their mean value is close to
\[
\frac{\sqrt{\ln\ln n\mathstrut}}{\ln n}.
\]

The paper ended with some theorems concerning the distribution of
upper and lower `peaks' of the graph of $S(t)$, i.e. the
distribution of quantities $S(\gamma_{n}+0)$, $S(\gamma_{n}-0)$
(theorems 11\,-14). As in the case of differences $\Delta_{n}$, the
major part of these peaks are of order $\sqrt{\ln\ln n\mathstrut}$.

\textsc{Notations.} Throughout the paper, $\vep$ denotes an
arbitrary small fixed positive number, $0<\vep<0.001$, $N\ge
N_{0}(\vep)>1$, $N$ is an integer, $M = \bigl[N^{27/82+\vep}\bigr]$,
$L = \ln\ln N$, $\theta, \theta_{1}, \theta_{2},\ldots$ are complex
numbers whose absolute value does not exceed one and which are,
generally speaking, different in different relations. All other
notations are explained in the text.

\textsc{Acknowledgments.} The author is grateful to Corresponding
Members of RAS, professors Yu.~V.~Nesterenko and A.~N.~Parshin for
comprehensive support and the attention to this paper, and to Dr.
T.~S.~Trudgian for sending his beautiful book \cite{Trudgian_2009}.

\pagebreak

\begin{center}
\textbf{\S 1. The Order of Growth of the Values}
$\boldsymbol{\Delta_{n}}$ \textbf{as} $\boldsymbol{n\to +\infty}$
\end{center}

In this part, some statements concerning the order of growth of
$\Delta_{n}$ are established. For the below, we need the following
definition.

\vspace{0.2cm}

\textbf{Definition 6.} Suppose $\gamma$ is an ordinate for $r$
different zeros of the Riemann zeta\,-function $\zeta(s)$ with
multiplicities $k_{1},\ldots, k_{r}$. Then the sum
\[
\kappa\,=\,k_{1}\,+\ldots +\,k_{r}
\]
is called the \emph{multiplicity of the ordinate} $\gamma$. In the
case $\gamma = \gamma_{l}$ we shall use the notation $\kappa_{l}$
instead $\kappa$. Thus in the case
\begin{equation}
\gamma_{l-1}<\gamma_{l} = \gamma_{l+1} = \ldots =
\gamma_{l+p-1}<\gamma_{l+p} \label{Lab9}
\end{equation}
we obviously have: $\kappa_{l} = \kappa_{l+1} = \ldots =
\kappa_{l+p-1} = p$. By basic properties of $S(t)$, if follows that
\begin{equation}
\kappa_{l}\,=\,S(\gamma_{l}+0)\,-\,S(\gamma_{l}-0).\label{Lab10}
\end{equation}

\textbf{Theorem 1.} $\Delta_{n}\,=\,O(\ln n)$ \emph{as} $n\to
+\infty$. \emph{If the Riemann hypothesis is true then}
\[
\Delta_{n}\,=\,O\biggl(\frac{\ln n}{\ln\ln n}\biggr).
\]

\emph{Proof.} Using the definition of $N(t)$, from (\ref{Lab1}) we
get
\begin{equation}
N(t_{m-1}-0)<N(\gamma_{n}+0)\le N(t_{m}+0).\label{Lab11}
\end{equation}
Now, $N(\gamma_{n}+0) = n+\theta_{n}(\kappa_{n}-1)$, where
$0\le\theta_{n}\le 1$. Indeed, the case $\kappa_{n} = 1$ is obvious.
Suppose the inequalities (\ref{Lab9}) hold for some $p\ge 2$. Then
for $n = l+s$, $0\le s\le p-1$ we have $\kappa_{n} = p$. Therefore,
\[
N(\gamma_{n}+0)\,=\,l+p-1\,=\,n+(p-s-1)\,=\,n+\kappa_{n}-s-1\,=\,n\,+\,\theta_{n}(\kappa_{n}-1),
\]
where
\[
0\,\le\,\theta_{n}\,=\,\frac{\kappa_{n}-s-1}{\kappa_{n}-1}\,\le\,1.
\]
Combining (\ref{Lab11}) with the Riemann\,-von Mangoldt formula we
conclude
\begin{multline*}
\frac{1}{\pi}\vth(t_{m-1})+1+S(t_{m-1}-0)\!<\!n+\theta_{n}(\kappa_{n}-1)\le
\frac{1}{\pi}\vth(t_{m})+1+S(t_{m}-0).
\end{multline*}
Using the definition of Gram points and replacing $\vth(t_{m-1})$,
$\vth(t_{m})$ by the quantities $\pi(m-2)$, $\pi(m-1)$, we get after
some transformations:
\[
-S(t_{m}+0)-\kappa_{n}\le m-n\le -S(t_{m-1}-0)+\kappa_{n}.
\]
Now the assertion follows from (\ref{Lab10}) and from classical
upper bounds for $|S(t)|$ (see, for example,
\cite{Karatsuba_Korolev_2005}; for the values of implied constants
in $O$'s, see the remark after lemma 7). This completes the proof.

\vspace{0.2cm}

The below theorem establishes the connection between the fractions
(\ref{Lab2}) and the quantities $\Delta_{n}$.

\textbf{Lemma 1.} \emph{As} $n\to +\infty$, \emph{the following
relation holds:}
\[
\Delta_{n}\,=\,q_{n}+\theta_{n}+O\biggl(\frac{\ln n}{n}\biggr),
\]
\emph{where} $0\le\theta_{n}\le 1$ \emph{and the constant in}
$O$\,\emph{-symbol is an absolute}.

\vspace{0.2cm}

\emph{Proof.} First, from (\ref{Lab1}) it follows that
\begin{equation}
\frac{t_{m-1}-t_{n}}{t_{n+1}-t_{n}}\,<\,q_{n}\,\le\,\frac{t_{m}-t_{n}}{t_{n+1}-t_{n}}.\label{Lab12}
\end{equation}
Next, the definition of Gram's points and the Lagrange's mean value
theorem imply the relation
\[
\pi\!\cdot\!(m-n) = \vth(t_{m})-\vth(t_{n}) =
\vth'(t_{n})(t_{m}-t_{n}) +
\tfrac{1}{2}\vth''(\xi)(t_{m}-t_{n})^{2},
\]
where $\xi$ lies between $t_{m}$ and $t_{n}$. Thus,
\begin{equation}
t_{m} -
t_{n}\,=\,\frac{\pi(m-n)}{\vth'(t_{n})}\cdot\frac{1}{1+r},\quad r =
\frac{\vth''(\xi)}{2\vth'(t_{n})}\,(t_{m}-t_{n}).\label{Lab13}
\end{equation}
Further, theorem 1 implies the following rough estimate
\[
t_{m}-t_{n}\ll t_{m}+t_{n}\ll \frac{n}{\ln n}.
\]
From the relations
\[
\vth''(\xi) = \frac{1}{2\xi} + O\biggl(\frac{1}{\xi^{3}}\biggr)\ll
\frac{1}{\xi}\ll \frac{1}{t_{n}}\ll \frac{\ln n}{n},\quad
\vth'(t_{n}) = \frac{1}{2}\ln\frac{t_{n}}{2\pi} +
O\biggl(\frac{1}{t_{n}}\biggr)\gg \ln n
\]
we obtain the inequality
\[
r\ll \frac{\ln n}{n}\cdot\frac{1}{\ln n}\cdot\frac{n}{\ln n}\ll
\frac{1}{\ln n}.
\]
Now the equality (\ref{Lab13}) implies more precise bounds for the
difference $t_{m}-t_{n}$ and for the quantity $r$, namely
\[
t_{m} - t_{n} =
\frac{\pi(m-n)}{\vth'(t_{n})}\,\biggl(1\,+\,O\biggl(\frac{1}{\ln
n}\biggr)\biggr)\ll \frac{|m-n|}{\ln n}\ll 1,\quad r\ll \frac{\ln
n}{n}\cdot\frac{1}{\ln n}\ll \frac{1}{n}.
\]
Thus we get
\begin{equation*}
t_{m} -
t_{n}\,=\,\frac{\pi(m-n)}{\vth'(t_{n})}\biggl(1\,+\,O\biggl(\frac{1}{n}\biggr)\biggr).
\end{equation*}
Further, the equalities
\[
t_{m-1}-t_{n} =
\frac{\pi(m-n-1)}{\vth'(t_{n})}\biggl(1\,+\,O\biggl(\frac{1}{n}\biggr)\biggr),
\]
\begin{equation}
t_{n+1} -
t_{n}\,=\,\frac{\pi}{\vth'(t_{n})}\biggl(1\,+\,O\biggl(\frac{1}{n\ln
n}\biggr)\biggr).\label{Lab14}
\end{equation}
are proved in the same way. Substituting all the expressions in
(\ref{Lab12}), we obtain
\begin{align*}
& q_{n}\,\le\,
\frac{\pi(m-n)}{\vth'(t_{n})}\cdot\frac{\vth'(t_{n})}{\pi}\biggl(1\,+\,O\biggl(\frac{1}{n}\biggr)\biggr)\,=\,\Delta_{n}\,+\,O\biggl(\frac{\ln
n}{n}\biggr),\\
& q_{n}\,>\,
\frac{\pi(m-n-1)}{\vth'(t_{n})}\cdot\frac{\vth'(t_{n})}{\pi}\biggl(1\,+\,O\biggl(\frac{1}{n}\biggr)\biggr)\,=\,\Delta_{n}-1+O\biggl(\frac{\ln
n}{n}\biggr).
\end{align*}
This completes the proof of the lemma.

\vspace{0.2cm}

\textbf{Theorem 2.} \emph{The inequalities}
\[
\max_{N<n\le N+M}\bigl(\pm\Delta_{n}\bigr)\,\ge\,c\biggl(\frac{\ln
N}{\ln\ln N}\biggr)^{\! 1/3},
\]
\emph{hold for some positive constant} $c = c(\vep)$.

\emph{Proof.} Let $Q = c_{1}\bigl(\frac{\ln N}{\ln\ln
N}\bigr)^{1/3}$ where the constant $c_{1}$ will be chosen later.
Suppose the inequalities
\begin{equation}
q_{n} = \frac{\gamma_{n}-t_{n}}{t_{n+1}-t_{n}} < Q.\label{Lab15}
\end{equation}
hold true for any $n$, $N<n\le N+M$. Taking $\tau_{n} =
t_{n}+Q(t_{n+1}-t_{n})$, we conclude from (\ref{Lab15}) that
$\gamma_{n}<\tau_{n}$. Then the Riemann\,-von Mangoldt formula
implies that the inequalities
\begin{equation}
n\le N(\gamma_{n}+0)\le N(\tau_{n}+0) =
\frac{1}{\pi}\vth(\tau_{n})+1+S(\tau_{n}+0)\label{Lab16}
\end{equation}
hold true for any $n$ from the interval under considering. Combining
the equality (\ref{Lab14}) with Lagrange's mean value theorem, we
have
\begin{multline*}
\vth(\tau_{n}) = \vth(t_{n}) + \vth'(t_{n})Q(t_{n+1}-t_{n}) +
\tfrac{1}{2}\vth''(\xi)Q^{2}(t_{n+1}-t_{n})^{2} =\\
=\,\vth(t_{n}) + \pi Q\biggl(1\,+\,O\biggl(\frac{1}{N\ln
N}\biggr)\biggr)\,+\,O\biggl(\frac{Q^{2}}{N\ln N}\biggr) =
\pi(n-1+Q) + O\bigl(N^{-1}\bigr)
\end{multline*}
for some $t_{n}<\xi< t_{n+1}$. Substituting this expression for
$\vth(\tau_{n})$ in (\ref{Lab16}), we obtain
\begin{equation}
S(\tau_{n}+0)\,\ge\,-Q\,+\,O\bigl(N^{-1}\bigr).\label{Lab17}
\end{equation}
Now let's show that the inequality (\ref{Lab17}) can't hold true for
every $n$ under condition $N<n\le N+M$. For arbitrary positive $x$,
by $t_{x}=t(x)$ we denote the unique solution of the equation
\begin{equation}
\vth(t_{x})\,=\,\pi\cdot(x-1)\label{Lab18}
\end{equation}
and set $\tau_{x} = t_{x} + Q(t_{x+1}-t_{x})$. Then
\[
\tau_{x}'\,=\,t_{x}'\,+\,Q(t_{x+1}'\,-\,t_{x}')\,=\,t_{x}'\,+\,Q\,t_{\eta}''
\]
for some $\eta$, $x<\eta<x+1$. By differentiating the equation
(\ref{Lab18}) with respect to $x$ twice we obtain
\begin{align*}
& t_{x}'\,=\,\frac{\pi}{\vth'(t_{x})},\\
&
t_{x}''\,=\,-\,\frac{(t_{x}')^{2}\vth''(t_{x})}{\vth'(t_{x})}\,=\,-\pi^{2}\,\frac{\vth''(t_{x})}{(\vth'(t_{x}))^{3}}\,\ll\,\frac{1}{x\ln^{2}x}.
\end{align*}
Therefore, the inequalities
\[
\tau_{x}'\,=\,\frac{\pi}{\vth'(t_{x})}\,-\,O\biggl(\frac{Q}{x\ln^{2}x}\biggr)\,=\,\frac{\pi}{\vth'(t_{x})}\biggl(1\,-\,O\biggl(\frac{Q}{x\ln
x}\biggr)\biggr)\,>\,\frac{\pi}{2\vth'(t_{x})}>0,
\]
holds for $N<x\le N+M$. Now it follows that the sequence $\tau_{n}$
increases monotonically in the interval under considering. Finally,
let's note that
\[
\tau_{n+1} - \tau_{n} =
\tau_{\zeta}'\,<\,\frac{\pi}{\vth'(t_{\zeta})}<\frac{\pi}{\vth'(t_{n})}
< \frac{3\pi}{\ln n}
\]
for some $n<\zeta < n+1$. Now let's show that if the point
$\tau_{n}$ is close to the minima of $S(t)$ then the value
$S(\tau_{n})$ is very large and negative. This will contradict to
(\ref{Lab17}).

From omega\,-theorem for $S(t)$ (see \cite{Karatsuba_Korolev_2005}),
it follows that there exists a point $\tau$ in the interval
$(\tau_{N},\tau_{N+M}]$ such that
\[
S(\tau)\,<\,-c_{0}\biggl(\frac{\ln N}{\ln\ln N}\biggr)^{1/3},
\]
for some constant $c_{0} = c_{0}(\vep) > 0$. Choosing $n$ from the
inequalities $\tau_{n}<\tau \le \tau_{n+1}$, we suppose that
$\gamma^{(1)}<\gamma^{(2)}<\ldots < \gamma^{(s)}$ are all different
ordinates of zeros of $\zeta(s)$ lying the interval
$(\tau_{n},\tau]$. Then, performing the increment of the function
$S(t)$ along this interval as the sum of increments along all
subintervals of continuity of $S(t)$ and the sum of jumps of $S(t)$
at the points of discontinuity $\gamma^{(1)},\gamma^{(2)},\ldots
,\gamma^{(s)}$, we get the following identity:
\begin{multline*}
S(\tau) -
S(\tau_{n}+0)\,=\,\bigl\{S(\gamma^{(1)}-0)-S(\tau_{n}+0)\bigr\}\,+\\
\bigl\{S(\gamma^{(1)}+0)-S(\gamma^{(1)}-0)\bigr\}\,+\,\bigl\{S(\gamma^{(2)}-0)-S(\gamma^{(1)}+0)\bigr\}\,+\,\ldots
\\ +
\bigl\{S(\gamma^{(s)}+0)-S(\gamma^{(s)}-0)\bigr\}\,+\,\bigl\{S(\tau)-S(\gamma^{(s)}+0)\bigr\}.
\end{multline*}
Obviously, this relation implies the following inequality:
\begin{multline}
S(\tau) -
S(\tau_{n}+0)\,\ge\,\bigl\{S(\gamma^{(1)}-0)-S(\tau_{n}+0)\bigr\}\,+\bigl\{S(\gamma^{(2)}-0)-S(\gamma^{(1)}+0)\bigr\}\,+\,\ldots
\\ +\,\bigl\{S(\tau)-S(\gamma^{(s)}+0)\bigr\}.\label{Lab19}
\end{multline}
Let $(a,b)$ be an interval of continuity of $S(t)$, and suppose that
$1<a<b$, $b-a<1$. Then
\begin{multline}
S(b) -
S(a)\,=\,(b-a)S'(\xi)\,=\,(b-a)\biggl(-\frac{1}{2}\ln\frac{\xi}{2\pi}\,+\,O\bigl(\xi^{-2}\bigr)\biggr)\,=\\
=\,(b-a)\biggl(-\frac{1}{2}\ln\frac{\xi}{2\pi}\,+\,O\bigl(a^{-2}\bigr)\biggr),\label{Lab20}
\end{multline}
for some $a<\xi<b$. Using the inequality
\[
\tau - \tau_{n}\le \tau_{n+1} - \tau_{n} < \frac{3\pi}{\ln{n}} < 1,
\]
from (\ref{Lab19}) and (\ref{Lab20}) we conclude
\begin{multline*}
S(\tau) -
S(\tau_{n}+0)\,\ge\,\bigl\{(\gamma^{(1)}-\tau_{n})+(\gamma^{(2)}-\gamma^{(1)})\,+\,\ldots\,+\,(\tau
-
\gamma^{(s)})\bigr\}\biggl(-\frac{1}{2}\ln\frac{\tau}{2\pi}\,+\,O\bigl(\tau^{-2}\bigr)\biggr)\,>\\
>\,-(\tau\,-\,\tau_{n})\ln{\tau}\,>\,-\frac{3\pi}{\ln
N}\cdot\ln{N}\,=\,-3\pi.
\end{multline*}
Hence
\[
S(\tau_{n}+0)\,\,\le\,S(\tau)+3\pi\,<\,-c_{0}\biggl(\frac{\ln
N}{\ln\ln N}\biggr)^{1/3} + 3\pi\,<\,-0.9c_{0}\biggl(\frac{\ln
N}{\ln\ln N}\biggr)^{1/3}.
\]
Now it is easy to see that the last inequality contradicts to
(\ref{Lab17}), if we choose $c_{1} = 0.8c_{0}$. Therefore, there
exists at least one number $n$ under condition $N<n\le N+M$, such
that the inequality (\ref{Lab15}) fails. By lemma 1, for such $n$ we
get
\[
\Delta_{n}\,>\,q_{n}-1+O\biggl(\frac{\ln{n}}{n}\biggr)\,\ge\,Q-1.1\,>\,c\biggl(\frac{\ln
N}{\ln\ln N}\biggr)^{1/3},
\]
where $c = 0.7c_{0}$.

The existence of large negative values of $\Delta_{n}$ is proved in
the same way. If we suppose that $q_{n}>-\,Q$ holds for every $n$
under considering then it follows that $\gamma_{n}>\sigma_{n}$ where
$\sigma_{n} = t_{n}-Q(t_{n+1}-t_{n})$, and hence
\[
N(\sigma_{n}+0)\,\le\,N(\gamma_{n}-0)\,<\,n.
\]
Transforming the left side as above we get
\begin{equation}
S(\sigma_{n}+0)\,<\,Q\,+\,O\bigl(N^{-1}\bigr).\label{Lab21}
\end{equation}

The application of omega\,-theorem for $S(t)$ yields the existence
of a point $\sigma$, $\sigma_{N}<\sigma\le \sigma_{N+M}$ such that
\[
S(\sigma)\,>\,c_{0}\biggl(\frac{\ln N}{\ln\ln N}\biggr)^{1/3}.
\]
Since the sequence $\sigma_{n}$ is monotonically increasing in the
interval under considering, it's possible to point out the number
$n$ such that $\sigma_{n-1}<\sigma\le \sigma_{n}$. Using the same
arguments as above and applying the inequality
$\sigma_{n}-\sigma_{n-1}<3\pi(\ln n)^{-1}$, we have
\begin{align*}
&
S(\sigma_{n}+0)-S(\sigma)\,>\,(\sigma_{n}-\sigma)\biggl(-\,\frac{1}{2}\ln\frac{\sigma_{n}}{2\pi}\,+\,O\bigl(\sigma_{n}^{-2}\bigr)\biggr)\,
>\,-(\sigma_{n}-\sigma)\ln{\sigma_{n}}\,>\,-3\pi,\\
& S(\sigma_{n}+0)\,>\,S(\sigma)-3\pi\,>\,0.9c_{0}\biggl(\frac{\ln
N}{\ln\ln N}\biggr)^{1/3}.
\end{align*}
The last relation contradicts to (\ref{Lab21}), if we choose $c_{1}
= 0.8c_{0}$. Therefore, there exists at least one number $n$,
$N<n\le N+M$, such that $q_{n}\le -Q$. For such $n$ we have
\[
\Delta_{n}\,\le\,q_{n}+1.1\,\le\,-\,Q\,+\,1.1\,<\,-\,c\biggl(\frac{\ln
N}{\ln\ln N}\biggr)^{1/3},
\]
where $c = 0.7c_{0}$. The theorem is completely proved.

\begin{center}
\textbf{\S 2. Selberg Hypothesis and the Moments of}
$\boldsymbol{\Delta_{n}}$
\end{center}

This section is devoted to the proof of Selberg hypothesis and to
the calculation of the moments of $\Delta_{n}$, that is, to the
calculation of the sums
\[
\sum\limits_{N<n\le N+M}\Delta_{n}^{a},\quad \sum\limits_{N<n\le
N+M}|\Delta_{n}|^{a}
\]
for different $a$. We need the following definition.

\textbf{Definition 7.} Suppose $a,b$ are an arbitrary real numbers,
$a<b$. Denote by $e(a,b)$ the number of solutions of the
inequalities $a\!<\!\Delta_{n}\le b$ with the condition $N\!<\! n\le
N+M$. Similarly, by $f(a,b)$ we denote the number of solutions of
the inequalities $a\!<\!\Delta(n)\le b$ under the same condition.

The below assertion plays the key role in what follows.

\textbf{Lemma 2.} \emph{The relation}
\[
e(a,b)\,=\,f(-(b+1),-(a+1))\,+\,\theta(|a|+|b|+2)
\]
\emph{holds true for any integers} $a$ \emph{and} $b$, $a<b$.

\emph{Proof.} First we get $e(a,b) = M-e_{1}-e_{2}$, where $e_{1}$
and $e_{2}$ denote, respectively, the numbers of solutions of
inequalities
\begin{equation}
\Delta_{n}\,\le\,a,\quad \Delta_{n}\,>\,b \label{Lab22}
\end{equation}
under the same condition $N\!<\!n\le N+M$. The first inequality in
(\ref{Lab22}) is equivalent to the following:
\begin{equation}
\gamma_{n}\,\le\,t_{n+a} \label{Lab23}
\end{equation}
Indeed, if $t_{m-1}\!<\!\gamma_{n}\le t_{m}$ and $\Delta_{n} =
m-n\le a$, then $m\le n+a$ and therefore $\gamma_{n}\le t_{m}\le
t_{n+a}$. Suppose now that (\ref{Lab23}) holds. Then for the number
$m$ defined above we have: $\gamma_{n}\le t_{m}$ $\le t_{n+a}$.
Hence, $m\le n+a$ and $\Delta_{n}\le a$.

By setting $\nu = n+a$ in (\ref{Lab23}), we obtain that $e_{1}$
equals to the number of solutions of the inequality
\begin{equation}
\gamma_{\nu - a}\,\le\,t_{\nu} \label{Lab24}
\end{equation}
with the condition $N+a<\nu\le N+M+a$. Hence, the difference between
$e_{1}$ and the number $f_{1}$ of solutions of (\ref{Lab24}) under
the condition $N<\nu\le N+M$ does not exceed $|a|$. Thus, $e_{1} =
f_{1}+\theta_{1}|a|$.

Using the same arguments we obtain that the inequalities
$\Delta_{n}>b$ and $\gamma_{n}>t_{n+b}$ are equivalent and the
quantity $e_{2}$ equals to the number of solutions $\nu$ such that
$N+b<\nu\le N+M+b$ and
\begin{equation}
\gamma_{\nu - b}\,>\,t_{\nu}. \label{Lab25}
\end{equation}
Hence, $e_{2} = f_{2}+\theta_{2}|b|$, where $f_{2}$ denotes the
number of solutions of (\ref{Lab25}) with the condition $N<\nu\le
N+M$.

Suppose the Gram point $t_{\nu}$ does not satisfy both (\ref{Lab24})
and (\ref{Lab25}). Then
\begin{equation}
\gamma_{\nu - b}\,\le\,t_{\nu}\,<\,\gamma_{\nu -a}. \label{Lab26}
\end{equation}
Let's show that (\ref{Lab26}) is equivalent to double inequality
\begin{equation}
-(b+1)<\Delta(\nu)\le -(a+1). \label{Lab27}
\end{equation}
Indeed, determining $m$ from the inequalities $\gamma_{m}\le
t_{\nu}<\gamma_{m+1}$, from (\ref{Lab26}) we conclude that $\nu
-b\le m\le \nu - a - 1$. Hence, (\ref{Lab27}) holds true. Next, for
$m$ defined above, from (\ref{Lab27}) it follows that $\nu - b \le
m\le \nu - a - 1$ and hence
\[
\gamma_{\nu - b}\le \gamma_{m}\le t_{\nu} < \gamma_{m+1}\le
\gamma_{\nu - a}.
\]

Therefore, the number of solutions of (\ref{Lab26}) satisfying the
condition $N<\nu\le N+M$, equals to $f(-(b+1),-(a+1))$. Thus,
\begin{multline*}
e(a,b)\,=\,M-e_{1}-e_{2}\,=\,M-f_{1}-f_{2}-\theta_{1}|a| -
\theta{2}|b|\,=\\
=\,f(-(b+1), -(a+1))\,+\,\theta(|a|+|b|).
\end{multline*}
Lemma is completely proved.

\vspace{0.2cm}

\textbf{Lemma 3.} \emph{For a real} $x$ \emph{let the quantity}
$\nu(x)$ \emph{denote the number of integers} $n$, $N<n\le N+M$,
\emph{satisfying the condition}
\[
\Delta(n)\,\le\,\frac{x}{\pi}\sqrt{\frac{L}{2}}.
\]
\emph{Then}
\[
\nu(x)\,=\,M\biggl(\frac{1}{\sqrt{2\pi}}\int_{-\infty}^{\,x}\!\!e^{-u^{2}/2}\,du\,+\,\theta\delta\biggr),
\]
\emph{where} $\delta = e^{22.4}\vep^{-1.5}(\ln L)^{-0.5}$.

\vspace{0.2cm}

For the proof, see \cite{Korolev_2010}.

\vspace{0.2cm}

\textbf{Theorem 3.} \emph{For any} $a$ \emph{and} $b$, $a<b$,
\emph{the number of solutions of the inequality} $a<\Delta_{n}\le b$
\emph{with the condition} $N<n\le N+M$ \emph{satisfies the relation}
\[
e(a,b)\,=\,M\biggl(\frac{1}{\sqrt{2\pi}}\int_{\alpha}^{\beta}\!e^{-u^{2}/2}\,du\,+\,\theta\Delta\biggr),
\]
\emph{where} $\alpha = \pi a\sqrt{2/L}$, $\beta = \pi b\sqrt{2/L}$,
$\Delta = 2.2e^{22.4}\vep^{-1.5}(\ln L)^{-0.5}$.

\vspace{0.2cm}

\emph{Proof.} Let $a$, $b$ be an integers and let $c$ be
sufficiently large constant such that the inequalities
$|\Delta_{n}|\le l$, $|\Delta(n)|\le l$ hold true for any $N<n\le
N+M$ with $l = [c\ln N]$. Then in the case $-l\le a < b \le l$ the
assertion follows from lemmas 2 and 3:
\begin{multline*}
e(a,b)\,=\,f(-(b+1),-(a+1))\,+\,2\theta_{1}l\,=\\
=\,M\biggl(\frac{1}{\sqrt{2\pi}}\int_{-\pi(b+1)\sqrt{2/L}}^{-\pi(a+1)\sqrt{2/L}}\!e^{-u^{2}/2}\,du\,+\,2\theta_{2}\delta\biggr)\,+\,2\theta_{1}l\,=\\
=\,M\biggl(\frac{1}{\sqrt{2\pi}}\int_{\pi a\sqrt{2/L}}^{\pi
b\sqrt{2/L}}\!e^{-u^{2}/2}\,du\,+\,\theta_{2}\sqrt{\pi/L}\,+\,2\theta_{2}\delta\biggr)\,+\,2\theta_{1}l\,=\\
=\,M\biggl(\frac{1}{\sqrt{2\pi}}\int_{\alpha}^{\beta}\!e^{-u^{2}/2}\,
du\,+\,2.1\theta\delta\biggr).
\end{multline*}
In the case $-l\le a < l <b$ the required statement follows from the
equality $e(a,b) = e(a,l)$ and from the estimate
\[
\frac{1}{\sqrt{2\pi}}\int_{\lambda}^{\beta}\!e^{-u^{2}/2}\,du\,<\,\lambda^{-1}e^{-\lambda^{2}/2}\,<\,0.01\delta,\quad
\lambda = \pi l\sqrt{2/L}.
\]
The cases $a<-l\le b \le l$, $a<-l<l<b$ are handle as above. If $a$
or $b$ is non\,-integer, then assertion follows from the relation
$e(a,b) = e([a], [b])$ and the above arguments. Theorem is
completely proved.

\vspace{0.2cm}

\textbf{Corollary.} \emph{Selberg hypothesis is true. Moreover},
\emph{if} $\Phi(x)>0$ \emph{is an arbitrary unbounded monotonically
increasing function, then the number of} $n$, $N<n\le N+M$,
\emph{that do not satisfy the condition}
\[
\frac{1}{\Phi(n)}\sqrt{\ln\ln n}<|\Delta_{n}|\le \Phi(n)\sqrt{\ln\ln
n},
\]
\emph{is of order not exceeding} $M(\Phi^{-1}(N)+\Delta)$.

\vspace{0.2cm}

Proof of this assertion repeats practically word\,-for\,-word the
proof of the corollary of theorem 4 in \cite{Korolev_2010}.

\vspace{0.2cm}

\textbf{Remark.} Theorem 3 asserts that the quantity $r(a,b)$ in the
relation
\[
e(a,b)\,=\,M\biggl(\frac{1}{\sqrt{2\pi}}\int_{\alpha}^{\beta}\!e^{-u^{2}/2}\,du\,+\,r(a,b)\biggr)
\]
obeys the estimate
\[
r(a,b)\,=\,O\bigl((\ln L)^{-0.5}\bigr)\,=\,O\bigl((\ln\ln\ln
N)^{-0.5}\bigr).
\]
for any $a$ and $b$. Suppose $\vep$ be any positive number, $\vep <
\tfrac{1}{2}$. Then for any $a$ we have
\[
0 = e(a+\vep,a+1-\vep) =
M\biggl(\frac{1}{\sqrt{2\pi}}\int_{\alpha_{1}}^{\alpha_{2}}\!e^{-u^{2}/2}\,du\,+\,r(a+\vep,a+1-\vep)\biggr),
\]
where $\alpha_{1} = \pi(a+\vep)\sqrt{2/L}$, $\alpha_{2} =
\pi(a+1-\vep)\sqrt{2/L}$. Hence,
\[
r(a+\vep,a+1-\vep)\,=\,-\frac{1}{\sqrt{2\pi}}\int_{\alpha_{1}}^{\alpha_{2}}\!e^{-u^{2}/2}du\,=\,O(\alpha_{2}-\alpha_{1})\,=\,O(L^{-0.5})\,=\,O\bigl((\ln\ln
N)^{-0.5}\bigr).
\]
Probably this estimate holds true for any $a$, $b$. It's interesting
to compare this assumption with one Selberg's assertion stated in
\cite[Theorem 2]{Selberg_1992}.

For the below, we need some new notations. For positive $a$, we set
\[
\vk(a)\,=\,\frac{2^{a}}{\sqrt{\pi}}\,\Gamma\biggl(\frac{a+1}{2}\biggr).
\]
Then for integer $k\ge 0$ we have:
\[
\vk(2k)\,=\,\frac{(2k)!}{k!},\quad
\vk(2k+1)\,=\,\frac{2^{2k+1}}{\sqrt{\pi}}\,k!.
\]
In what follows, we shall use the obvious estimations
\[
\vk(2k)\,<\,\frac{3}{2}\biggl(\frac{4k}{e}\biggr)^{k},\quad
\frac{1}{\vk(2k)}\,<\,\biggl(\frac{e}{4k}\biggr)^{k},
\]
without special comments, and the same for inequalities
\[
k\,\le\,\bigl(\!\sqrt[3\;]{3}\bigr)^{k},\quad
k\sqrt{k}\,\le\,\bigl(\!\sqrt{3}\bigr)^{k},\quad
k^{2}\,\le\,\bigl(\!\sqrt[3\;]{9}\bigr)^{k}.
\]
Finally, we set $A = e^{21}\vep^{-1.5}$, $B = A^{2}e^{-8} =
e^{34}\vep^{-3}$. In order to prove Selberg's formulae, we need some
auxiliary assertions from \cite{Korolev_2010}.

\vspace{0.2cm}

\textbf{Lemma 4.} \emph{Let} $k$ \emph{be an integer with the
condition} $1\le k\le \sqrt{L}$. \emph{Then the following relation
holds}:
\[
\sum\limits_{N<n\le
N+M}\Delta^{2k}(n)\,=\,\frac{\vk(2k)}{(2\pi)^{2k}}\,ML^{k}\bigl(1\,+\,\theta
A^{k}L^{-0.5}\bigr).
\]
\emph{Moreover}, \emph{we have the estimate}
\[
\biggl|\sum\limits_{N<n\le
N+M}\Delta^{2k-1}(n)\biggr|\,\le\,\frac{3.5}{\sqrt{B}}\,(Bk)^{k}ML^{k-1}.
\]

\vspace{0.2cm}

\textbf{Lemma 5.} \emph{Let} $a$ \emph{be an arbitrary number with
the condition}
\[
0<a\le\frac{e\ln L}{10A\ln\ln L}.
\]
\emph{Then the following asymptotic formula holds}:
\[
\sum\limits_{N<n\le
N+M}|\Delta(n)|^{a}\,=\,\frac{\vk(a)}{(2\pi)^{a}}\,ML^{0.5a}\bigl(1\,+\,\theta
c(A^{-1}\ln L)^{-\gamma}\bigr).
\]
\emph{Furthermore, in the case} $0<a\le 1$ \emph{the quantities} $c$
\emph{and} $\gamma$ \emph{can be set to be equal to} $90$ \emph{and}
$0.5a$, \emph{respectively, and in the case} $a>1$ \emph{to}
$2^{15.4}$ \emph{and} $\tfrac{1}{2}+\bigl\{\frac{a-1}{2}\bigr\}$.

Similarly to the assertions of lemmas 4 and 5, the below analogues
of Selberg's formulae are uniform to the parameter $k$.

\textbf{Theorem 4.} \emph{Let} $k$ \emph{be an integer with the
condition} $1\le k\le \sqrt{L}$. \emph{Then the following equality
holds}:
\[
\sum\limits_{N<n\le
N+M}\Delta_{n}^{2k}\,=\,\frac{\vk(2k)}{(2\pi)^{2k}}\,ML^{k}\bigl(1\,+\,1.1\theta
A^{k}L^{-0.5}\bigr).
\]
\emph{In particular}, \emph{for fixed} $k$ \emph{we get}
\[
\sum\limits_{N<n\le
N+M}\Delta_{n}^{2k}\,=\,\frac{(2k)!}{k!}\,\frac{M}{(2\pi)^{2k}}\,(\ln\ln
N)^{k}\,+\,O\bigl(M(\ln\ln N)^{k - 1/2}\bigr).
\]

\vspace{0.2cm}

\textbf{Corollary.} \emph{Under the same restrictions on} $k$,
\emph{the following estimate holds}:
\[
\sum\limits_{N<n\le
N+M}\Delta_{n}^{2k}\,\le\,\frac{1.1}{A}\,\frac{\vk(2k)}{(2\pi)^{2k}}\,M(AL)^{k}.
\]

\textbf{Theorem 5.} \emph{Let} $k$ \emph{be an integer with the
condition} $1\le k\le \sqrt{L}$. \emph{Then the following estimate
holds}:
\[
\biggl|\sum\limits_{N<n\le
N+M}\Delta_{n}^{2k-1}\biggr|\,\le\,\frac{3.6}{\sqrt{B}}\,(Bk)^{k}ML^{k-1}.
\]
\emph{In particular}, \emph{for fixed} $k$ \emph{we have}:
\[
\sum\limits_{N<n\le N+M}\Delta_{n}^{2k-1}\,=\,O\bigl(M(\ln\ln
N)^{k-1}\bigr).
\]

\vspace{0.2cm}

\textbf{Theorem 6.} \emph{Let} $a$ \emph{be an arbitrary number with
the condition}
\[
0<a\le\frac{e\ln L}{10A\ln\ln L}.
\]
\emph{Then the following asymptotic formula holds}:
\[
\sum\limits_{N<n\le
N+M}|\Delta_{n}|^{a}\,=\,\frac{\vk(a)}{(2\pi)^{a}}\,ML^{0.5a}\bigl(1\,+\,\theta
c(A^{-1}\ln L)^{-\gamma}\bigr).
\]
\emph{Furthermore, in the case} $0<a\le 1$ \emph{the quantities} $c$
\emph{and} $\gamma$ \emph{can be set to be equal to} $91$ \emph{and}
$0.5a$, \emph{respectively, and in the case} $a>1$ \emph{to}
$2^{15.5}$ \emph{and} $\tfrac{1}{2}+\bigl\{\frac{a-1}{2}\bigr\}$.
\emph{In particular}, \emph{in the case} $a = 1$ \emph{we get the
equality}:
\[
\sum\limits_{N<n\le
N+M}|\Delta_{n}|\,=\,\frac{M}{\pi\sqrt{\pi}}\sqrt{\ln\ln
N}\bigl(1\,+\,O((\ln\ln\ln N)^{-1/2})\bigr).
\]

\emph{Proof.} Theorems 4\,-6 can be proved in a similar way. Using
the fact that the numbers $\Delta_{n}$ are integral, for any $a>0$
we obtain
\[
\sum\limits_{N<n\le N+M}|\Delta_{n}|^{a}\,=\,\sum\limits_{\nu =
1}^{l}\nu^{a}\bigl(e(\nu -1,\nu)\,+\,e(-(\nu+1),-\nu)\bigr).
\]
The number $l = [c\ln N]$ is chosen such that the absolute values of
the quantities $\Delta_{n}$, $\Delta(n)$ do not exceed $l$ for
$N<n\le N+M$. By lemma 2, for $\nu\ge 2$ we get
\[
e(\nu -1,\nu)\,=\,f(-(\nu+1),-\nu)\,+\,\theta_{1}(2\nu-1),\quad
e(-(\nu+1),-\nu))\,=\,f(\nu-1,\nu)\,+\,\theta_{2}(2\nu+1).
\]
Hence,
\begin{multline*}
\sum\limits_{N<n\le N+M}|\Delta_{n}|^{a}\,=\,\sum\limits_{\nu =
1}^{l}\nu^{a}\bigl(f(\nu-1,\nu)\,+\,f(-(\nu+1),-\nu)\,+\,4\theta_{3}\nu\bigr)\,=\\
=\,\sum\limits_{N<n\le
N+M}|\Delta(n)|^{a}\,+\,4\theta_{4}\sum\limits_{\nu =
1}^{l}\nu^{a+1}.
\end{multline*}
Now theorems 4 and 6 follow from lemmas 4 and 5, respectively, and
from obvious bound
\[
4\sum\limits_{\nu =
1}^{l}\nu^{a+1}\,<\,\frac{4(l+1)^{a+2}}{a+2}\,<\,\bigl(1.5c\ln
N\bigr)^{a+2}.
\]
Further, for integer $k\ge 1$ we have:
\begin{multline*}
\sum\limits_{N<n\le N+M}\!\!\Delta_{n}^{2k-1}\,=\,\sum\limits_{\nu =
1}^{l}\bigl(\nu^{2k-1}e(\nu-1,\nu)\,+\,(-\nu)^{2k-1}e(-(\nu+1),-\nu)\bigr)\,=\\
=\,\sum\limits_{\nu =
1}^{l}\nu^{2k-1}\bigl(e(\nu-1,\nu)\,-\,e(-(\nu+1),-\nu)\bigr)\,=\\
=\,\sum\limits_{\nu =
1}^{l}\nu^{2k-1}\bigl(f(-(\nu+1),-\nu)\,-\,f(\nu-1,\nu)\,+\,4\theta_{1}\nu\bigr)\,=\\
=\,-\,\sum\limits_{\nu =
1}^{l}\bigl(\nu^{2k-1}f(\nu-1,\nu)\,+\,(-\nu)^{2k-1}f(-(\nu+1),-\nu)\bigr)\,+\,\theta_{2}\bigl(1.5c\ln
N\bigr)^{a+2}\,=\\
=\,-\sum\limits_{N<n\le
N+M}\!\!\Delta^{2k-1}(n)\,+\,\theta_{3}\bigl(1.5c\ln N\bigr)^{a+2}.
\end{multline*}
Thus, the assertion of Theorem 5 is a corollary of lemma 4.

\begin{center}
\textbf{\S 3. The Distribution of Differences}
$\boldsymbol{\gamma_{n}-t_{n}}$
\end{center}

Lemma 1 implies that the difference between quantities $\Delta_{n}$
and $q_{n} = \frac{\gamma_{n}-t_{n}}{t_{n+1}-t_{n}}$ does not exceed
$O(1)$. This fact together with the assertions of theorems 3\,-6
allow us to find approximately the distribution function for the
differences $t_{n}-\gamma_{n}$ and to calculate the moments of these
quantities.

\textbf{Lemma 6.} \emph{For any} $n$ \emph{with the condition}
$N<n\le N+M$ \emph{the following equality holds}
\[
\gamma_{n}-t_{n}\,=\,\frac{\pi}{\vth'(t_{N})}\bigl(\Delta_{n}\,+\,\vep_{n}\bigr),
\]
\emph{where} $|\vep_{n}|\le 1.01$.

\emph{Proof.} Using (\ref{Lab14}) and the relation
\[
\frac{1}{\vth'(t_{n})}\,=\,\frac{\pi}{\vth'(t_{N})}\biggl(1\,+\,O\biggl(\frac{M}{N\ln
N}\biggr)\biggr),
\]
we get:
\[
\gamma_{n}-t_{n}\,=\,q_{n}(t_{n+1}-t_{n})\,=\,\frac{\pi
q_{n}}{\vth'(t_{n})}\biggl(1\,+\,O\biggl(\frac{1}{N\ln
N}\biggr)\!\biggr)\,=\,\frac{\pi
q_{n}}{\vth'(t_{N})}\biggl(1\,+\,O\biggl(\frac{M}{N\ln
N}\biggr)\!\biggr).
\]
Expressing $q_{n}$ by $\Delta_{n}$, from lemma 1 we obtain
\[
\gamma_{n}-t_{n}\,=\,\frac{\pi}{\vth'(t_{N})}\biggl(\Delta_{n}-\theta_{n}+O\biggl(\frac{\ln
N}{N}\biggr)\!\biggr)\biggl(1\,+\,O\biggl(\frac{M}{N\ln
N}\biggr)\!\biggr)\,=\,\frac{\pi}{\vth'(t_{N})}\bigl(\Delta_{n}\,+\,\vep_{n}\bigr),
\]
where
\[
|\vep_{n}|\,=\,\bigl|\theta_{n}\,+\,O(MN^{-1})\bigr|\,\le\,1.01.
\]
The lemma is proved.

\vspace{0.2cm}

\textbf{Theorem 7.} \emph{Let} $k$ \emph{be an integer with the
condition} $1\le k\le \sqrt{L}$. \emph{Then the following equality
holds}
\[
\sum\limits_{N<n\le
N+M}(\gamma_{n}-t_{n})^{2k}\,=\,\frac{\vk(2k)}{(2\vth'(t_{N}))^{2k}}\,ML^{k}\bigl(1\,+\,0.4\theta(4A\sqrt[6\;]{3})^{k}L^{-0.5}\bigr).
\]

\vspace{0.2cm}

\textbf{Proof.} Using the equality
\[
(a+b)^{2k}\,=\,a^{2k}\,+\,\theta\,
k2^{2k-1}\bigl(|a|^{2k-1}|b|\,+\,|b|^{2k}\bigr)
\]
and defining the quantities $\vep_{n}$ as in lemma 6, we find
\begin{align*}
& (\Delta_{n}+\vep_{n})^{2k}\,=\,\Delta_{n}^{2k}\,+\,\theta\,k2^{2k-1}\bigl(|\Delta_{n}|^{2k-1}|\vep_{n}|\,+\,|\vep_{n}|^{2k}\bigr)\,=\\
&
=\,\Delta_{n}^{2k}\,+\,\theta\,k2^{2k}\bigl(|\Delta_{n}|^{2k-1}\,+\,(1.01)^{2k}\bigr),\\
& \sum\limits_{N<n\le
N+M}(\Delta_{n}+\vep_{n})^{2k}\,=\,\Sigma_{1}\,+\,\theta\,k
2^{2k}\bigl(\Sigma_{2}\,+\,(1.01)^{2k}M\bigr),
\end{align*}
where
\[
\Sigma_{1}\,=\,\sum\limits_{N<n\le N+M}\Delta_{n}^{2k},\quad
\Sigma_{2}\,=\,\sum\limits_{N<n\le N+M}|\Delta_{n}|^{2k-1}.
\]
Further,
\[
M\,=\,\frac{\vk(2k)}{(2\pi)^{2k}}\,ML^{k}\delta_{1},\quad
\delta_{1}\,=\,\frac{(2\pi)^{2k}}{\vk(2k)}\,\frac{1}{L^{k}}\,<\,\biggl(\frac{4\pi^{2}}{L}\biggr)^{k}\biggl(\frac{e}{4k}\biggr)^{k}\,=\,\biggl(\frac{\pi^{2}e}{kL}\biggr)^{k}
\,\le\,\frac{\pi^{2}e}{L},
\]
and, furthermore,
\[
\delta_{1}^{1/(2k)}\,<\,\sqrt{\frac{\pi^{2}e}{kL}}.
\]
Applying H\"{o}lder's inequality together with the corollary of
theorem 4 we get
\begin{align*}
& k2^{2k}\Sigma_{2}\,\le\,
k2^{2k}M^{1/(2k)}\Sigma_{1}^{1-1/(2k)}\,=\,\frac{\vk(2k)}{(2\pi)^{2k}}\,ML^{k}\delta_{2},\\
&
\delta_{2}\,=\,k2^{2k}\sqrt{\frac{\pi^{2}e}{kL}}\,\frac{1.1A^{k}}{\sqrt{A}}\,<\,\frac{0.1(4A\sqrt[6\;]{3})^{k}}{\sqrt{L}}.
\end{align*}
Returning to the initial sum, we obtain
\begin{align*}
& \sum\limits_{N<n\le
N+M}(\Delta_{n}+\vep_{n})^{2k}\,=\,\frac{\vk(2k)}{(2\pi)^{2k}}\,ML^{k}(1\,+\,\theta\delta),\\
&
\delta\,=\,\frac{1.1A^{k}}{\sqrt{L}}\,+\,\frac{0.1(4A\sqrt[6\;]{3})^{k}}{\sqrt{L}}\,+\,k(2.02)^{2k}\,\frac{\pi^{2}e}{L}\,<\,\frac{0.4(4A\sqrt[6\;]{3})^{k}}{\sqrt{L}},\\
& \sum\limits_{N<n\le
N+M}\!\!(\gamma_{n}-t_{n})^{2k}=\biggl(\frac{\pi}{\vth'(t_{N})}\biggr)^{2k}\!\!\!\sum\limits_{N<n\le
N+M}(\Delta_{n}+\vep_{n})^{2k}=\frac{\vk(2k)}{(2\vth'(t_{N}))^{2k}}\,ML^{k}(1+\theta
\delta).
\end{align*}
The theorem is proved.

\vspace{0.2cm}

\textbf{Remark.} If we set
$e_{n}\,=\,(\gamma_{n}-t_{n})\vth'(t_{N})\sqrt{2/L}$, then the
assertion of theorem 7 can be represented in the form:
\[
\sum\limits_{N<n\le
N+M}e_{n}^{2k}\,=\,\frac{(2k)!}{k!}\,\frac{M}{2^{k}}\,\bigl(1\,+\,0.4\theta(4A\sqrt[6\;]{3})L^{-0.5}\bigr).
\]

\vspace{0.2cm}

\textbf{Theorem 8.} \emph{Let} $k$ \emph{be an integer with the
condition} $1\le k\le \sqrt{L}$. \emph{Then the following estimation
holds}:
\[
\biggl|\sum\limits_{N<n\le
N+M}(\gamma_{n}\,-\,t_{n})^{2k-1}\biggr|\,<\,\frac{3.7}{\sqrt{B}}\biggl(\frac{\pi}{\vth'(t_{N})}\biggr)^{2k-1}(Bk)^{k}ML^{k-1}.
\]

\vspace{0.2cm}

\emph{Proof.} By the relation
\[
(a-b)^{2k-1}\,=\,a^{2k-1}\,+\,\theta\,k2^{2k-1}\bigl(a^{2k-2}|b|\,+\,|b|^{2k-1}\bigr),
\]
we obtain:
\[
\sum\limits_{N<n\le
N+M}(\Delta_{n}+\vep_{n})^{2k-1}=\!\!\sum\limits_{N<n\le
N+M}\!\!\Delta_{n}^{2k-1}\,+\,\theta_{1}\bigl(k\,
2^{2k}\!\!\sum\limits_{N<n\le
N+M}\Delta_{n}^{2k-2}+k(2.02)^{2k}M\bigr).
\]
First, the corollary of theorem 4 implies the estimations
\begin{align*}
& k\,2^{2k}\!\!\sum\limits_{N<n\le
N+M}\Delta_{n}^{2k-2}\,\le\,k\,2^{2k}\,\frac{1.1}{A}\,\frac{\vk(2k-2)}{(2\pi)^{2k-2}}M(AL)^{k-1}\,=\,
\frac{3.6}{\sqrt{B}}\,(Bk)^{k}ML^{k-1}\delta_{1},
\end{align*}
where
\begin{multline*}
\delta_{1}\,=\,\frac{\sqrt{B}}{3.6}\cdot\frac{1.1k\,2^{2k}}{(Bk)^{k}}\,\frac{(4A)^{k}}{A^{2}(2\pi)^{2k-2}}\,\vk(2k-2)\,=\,\frac{1.1}{3.6}\,\frac{\sqrt{B}}{A^{2}}\cdot
4\pi^{2}
\biggl(\frac{4A}{4\pi^{2}Bk}\biggr)^{k}\frac{k^{2}\vk(2k-2)}{2k(2k-1)}\,<\\
<\,\frac{1.1}{3.6}\,\frac{\sqrt{B}}{A^{2}}\cdot
2\pi^{2}\biggl(\frac{A}{\pi^{2}Bk}\biggr)^{k}\,\frac{3}{2}\biggl(\frac{4k}{e}\biggr)^{k}\,=\,\frac{1.1}{1.2}\,\pi^{2}\,\frac{\sqrt{B}}{A^{2}}
\biggl(\frac{4A}{\pi^{2}Be}\biggr)^{k}\,<\,\pi^{2}\,\frac{\sqrt{B}}{A^{2}}\,\frac{4A}{\pi^{2}Be}\,=\,\frac{4e^{-5}}{B}.
\end{multline*}
Further,
\begin{align*}
&
k(2.02)^{2k}M\,=\,\frac{3.6}{\sqrt{B}}(Bk)^{k}ML^{k-1}\delta_{2},\\
&
\delta_{2}\,=\,\frac{\sqrt{B}}{3.6}\cdot\frac{k}{(Bk)^{k}}\,\frac{(2.02)^{2k}}{L^{k-1}}\,<\,\frac{\sqrt{B}}{3.6}\,\frac{1}{(kL)^{k-1}}\biggl(\frac{4.1}{B}\biggr)^{k}\,\le\,
\frac{4.1}{3.6}\,\frac{1}{\sqrt{B}}\,<\,\frac{1.2}{\sqrt{B}}.
\end{align*}
Finally, combining all these estimations and using the theorem 5, we
get
\begin{align*}
& \biggl|\sum\limits_{N<n\le
N+M}(\Delta_{n}+\vep_{n})^{2k-1}\biggr|\,<\,\frac{3.6}{\sqrt{B}}(Bk)^{k}ML^{k-1}\biggl(1\,+\,\frac{4e^{-5}}{B}\,+\,\frac{1.2}{\sqrt{B}}\biggr)\,<\\
& <\,\frac{3.7}{\sqrt{B}}(Bk)^{k}ML^{k-1}.
\end{align*}
The last relation implies the required assertion.

\vspace{0.2cm}

\textbf{Remark.} The estimate of theorem 8 can be represented in the
form:
\[
\biggl|\sum\limits_{N<n\le
N+M}e_{n}^{2k-1}\biggr|\,<\,\frac{2.4}{\sqrt{B}}\,(Be\pi^{2})^{k}\cdot
\frac{(2k-1)!}{k!}\,ML^{-0.5}.
\]

\vspace{0.2cm}

\textbf{Theorem 9.} \emph{Let} $a$ \emph{be an arbitrary number with
the condition}
\[
0<a<\frac{e\ln{L}}{10A\ln\ln{L}}.
\]
\emph{Then the following asymptotic formula holds}:
\[
\sum\limits_{N<n\le
N+M}|\gamma_{n}-t_{n}|^{a}\,=\,\frac{\vk(a)}{(2\vth'(t_{N}))^{a}}ML^{0.5}\bigl(1\,+\,\theta
c_{1}(A^{-1}\ln{L})^{-\gamma}\bigr).
\]
\emph{Furthermore, in the case} $0<a \le 1$ \emph{the quantities}
$c_{1}$ \emph{and} $\gamma$ \emph{can be set to be equal} $97$
\emph{and} $0.5a$, \emph{respectively}, \emph{and in the case} $a>
1$ \emph{to} $2^{15.6}$ \emph{and}
$\tfrac{1}{2}+\bigl\{\frac{a-1}{2}\bigr\}$. \emph{In particular},
\emph{in the case} $a = 1$ \emph{we get the equality}:
\[
\sum\limits_{N<n\le
N+M}|\gamma_{n}\,-\,t_{n}|\,=\,\frac{2M}{\sqrt{\pi}}\,\frac{\sqrt{\ln\ln
N}}{\ln N}\bigl(1\,+\,O((\ln\ln\ln N)^{-0.5})\bigr).
\]

\vspace{0.2cm}

\emph{Proof.} First consider the case $0<a\le 1$. Let's show that
\begin{equation}
|\Delta_{n}\,+\,\vep_{n}|^{a}\,=\,|\Delta_{n}|^{a}\,+\,1.01\theta.\label{Lab28}
\end{equation}
Indeed, in the case $\Delta_{n} = 0$ the equality (\ref{Lab28})
follows from the estimate $|\vep_{n}|\le 1.01$. If $\Delta_{n} = 1$
then $|\Delta_{n}+\vep_{n}|^{a} = 2.01\theta$ and hence
\[
|\Delta_{n}+\vep_{n}|^{a}\,=\,1+2.01\theta-1\,=\,1+1.01\theta_{1}\,=\,|\Delta_{n}|^{a}\,+\,1.01\theta_{1}.
\]
If $\Delta_{n}\ge 2$ then $\Delta_{n}+\vep_{n}>0$. By Lagrange's
mean value theorem,
\[
|\Delta_{n}+\vep_{n}|^{a}\,=\,(\Delta_{n}+\vep_{n})^{a}\,=\,\Delta_{n}^{a}\,+\,a\vep_{n}(\Delta_{n}+\theta\vep_{n})^{a-1}.
\]
Since $|\Delta_{n}+\theta\vep_{n}|\ge 2-1.01 = 0.99$, we get:
\begin{align*}
&
|a\vep_{n}(\Delta_{n}+\theta\vep_{n})^{a-1}|\,\le\,1.01a\cdot(0.99)^{a-1}\,\le\,1.01,\\
& |\Delta_{n}+\vep_{n}|^{a}\,=\,|\Delta_{n}|^{a}\,+\,1.01\theta.
\end{align*}
The case of negative $\Delta_{n}$ is handled as above. Therefore,
\[
\sum\limits_{N<n\le
N+M}|\Delta_{n}+\vep_{n}|^{a}\,=\,\sum\limits_{N<n\le
N+M}|\Delta_{n}|^{a}\,+\,1.01\theta M.
\]
Now the required assertion follows from theorem 6, from the equality
\[
1.01M\,=\,\frac{\vk(a)}{(2\pi)^{a}}\,ML^{0.5a}\delta_{1}
\]
and from the estimations
\[
\delta_{1}\,=\,\frac{\pi^{a+0.5}}{\Gamma\Bigl(\D\frac{a+1}{2}\Bigr)}\,\frac{1.01}{L^{0.5a}}\,\le\,1.01\pi\sqrt{\pi}L^{-\,\gamma}\,<\,5.7L^{-\,\gamma}\,<\,
6\bigl(A^{-1}\ln{L}\bigr)^{-\,\gamma}.
\]
Suppose now $a>1$. Using the same arguments as above, we get
\[
|\Delta_{n}+\vep_{n}|^{a}\,=\,|\Delta_{n}|^{a}\,+\,\theta\,(0.6a2^{a}|\Delta_{n}|^{a-1}\,+\,a(2.02)^{a}).
\]
Summing the both parts over $n$, we obtain:
\[
\sum\limits_{N<n\le
N+M}\!\!|\Delta_{n}+\vep_{n}|^{a}=\sum\limits_{N<n\le
N+M}\!\!|\Delta_{n}|^{a}\,+\,\theta\bigl(0.6a\,2^{a}\sum\limits_{N<n\le
N+M}|\Delta_{n}|^{a-1}\,+\,a(2.02)^{a}M\bigr).
\]
Further, theorem 6 implies that the last sum over $n$ does not
exceed
\[
\frac{\Gamma(0.5a)}{\pi^{a-0.5}}\,ML^{(a-1)/2}d,\quad
d\,=\,1\,+\,c_{1}\bigl(A^{-1}\ln{L}\bigr)^{-\,\gamma},
\]
where $c_{1} = 91$, $\gamma = \tfrac{1}{2}(a-1)$ in the case $1<a\le
2$, and $c_{1} = 2^{15.5}$, $\gamma =
\tfrac{1}{2}+\bigl\{\frac{a-2}{2}\bigr\} =
\tfrac{1}{2}+\bigl\{\frac{a}{2}\bigr\}$ in the case $a>2$. It's easy
to note that the quantity $d$ does not exceed $92$ for any $a>1$.
Hence,
\[
0.6a\,2^{a}\sum\limits_{N<n\le N+M}|\Delta_{n}|^{a-1}\,\le\, 92\cdot
0.6a\,2^{a}\frac{\Gamma(0.5a)}{\pi^{a-0.5}}\,ML^{(a-1)/2}\,=\,\frac{\vk(a)}{(2\pi)^{a}}\,ML^{0.5a}\delta_{1},
\]
where
\begin{multline*}
\delta_{1}\,=\,\frac{\pi^{a+0.5}}{\Gamma(0.5(a+1))}\,\frac{\Gamma(0.5a)}{\pi^{a-0.5}}\,92\cdot
0.6a\,2^{a}L^{-0.5}\,<\,348\,\frac{\Gamma(0.5a+1)}{\Gamma(0.5(a+1))}\,2^{a}L^{-0.5}\,<\\
<\,348(a+1)2^{a}L^{-0.5}\,<\,L^{-0.49}.
\end{multline*}
Finally, we have:
\[
0.5a(2.02)^{a}M\,=\,\frac{\vk(a)}{(2\pi)^{a}}\,ML^{0.5a}\delta_{2},
\]
\[
\delta_{2}\,=\,\frac{a(2.02)^{a}\pi^{a+0.5}}{\Gamma(0.5(a+1))L^{0.5a}}\,\le\,2a\sqrt{\pi}\biggl(\frac{2.02\pi}{\sqrt{L}}\biggr)^{a}\,
<\,\frac{4.04\pi\sqrt{\pi}}{\sqrt{L}}\,<\,\frac{23}{\sqrt{L}}.
\]
Thus,
\[
\sum\limits_{N<n\le
N+M}|\Delta_{n}+\vep_{n}|^{a}\,=\,\frac{\vk(a)}{(2\pi)^{a}}\,ML^{0.5a}(1\,+\,\theta\delta),
\]
where
\begin{align*}
&
|\delta|\,\le\,2^{15.5}\bigl(A^{-1}\ln{L}\bigr)^{-\,\gamma}\,+\,L^{-0.49}\,+\,23L^{-0.5}\,<\,2^{15.6}\bigl(A^{-1}\ln{L}\bigr)^{-\,\gamma},\\
& \sum\limits_{N<n\le
N+M}|\gamma_{n}-t_{n}|^{a}\,=\,\frac{\vk(a)}{(2\vth'(t_{N}))^{a}}\,ML^{0.5a}(1\,+\,\theta\delta).
\end{align*}
The theorem is proved.

\vspace{0.2cm}

The assertions of theorems 7 and 8 allow one to find approximately
the distribution function for the differences $\gamma_{n}-t_{n}$.

\vspace{0.2cm}

\textbf{Theorem 10.} \emph{The quantity} $\nu(x)$ \emph{of the
numbers} $n$, $N<n\le N+M$ \emph{that satisfy the condition}
\[
e_{n}\,=\,(\gamma_{n}-t_{n})\vth'(t_{N})\sqrt{2/L}\,\le\,x
\]
\emph{obeys the following relation}:
\[
\nu(x)\,=\,M\biggl(\frac{1}{\sqrt{2\pi}}\int_{-\infty}^{x}e^{-u^{2}/2}\,du\,+\,\frac{\theta
A}{\sqrt{\ln{L}}}\biggr).
\]

\vspace{0.2cm}

\emph{Proof.} We give only a sketch, since the proof of this
assertion repeats practically word\,-for\,-word the proof of theorem
4 in \cite{Korolev_2010} (lemma 3 in present paper). That's why the
most part of calculations is missed here.

Suppose $f(t)$ be a characteristic function of the discrete random
quantity with the values $e_{n}$, $N<n\le N+M$, that is, the sum
\[
f(t)\,=\,\frac{1}{M}\sum\limits_{N<n\le N+M}\exp{(ite_{n})}.
\]
Taking an integer $K>1$ whose precise value is chosen below, by
corollaries of theorems 7 and 8 we get
\begin{multline*}
f(t)\,=\,\frac{1}{M}\sum\limits_{N<n\le
N+M}\bigl(\cos{(te_{n})}\,+\,i\sin{(te_{n})}\bigr)\,=\\
=\,\frac{1}{M}\sum\limits_{N<n\le N+M}\biggl(\sum\limits_{k =
0}^{K-1}\frac{(-1)^{k}}{(2k)!}\,(te_{n})^{2k}\,+\,i\sum\limits_{k =
1}^{K}\frac{(-1)^{k-1}}{(2k-1)!}(te_{n})^{2k-1}\,+\,2\theta\,\frac{(te_{n})^{2K}}{(2K)!}\biggr)\,=\\
=\,\sum\limits_{k =
0}^{K-1}\frac{(-1)^{k}t^{2k}}{(2k)!}\,\frac{1}{M}\sum\limits_{N<n\le
N+M}\!\!e_{n}^{2k}\,+\,i\sum\limits_{k =
1}^{K}\frac{(-1)^{k-1}t^{2k-1}}{(2k-1)!}\sum\limits_{N<n\le
N+M}\!\!e_{n}^{2k-1}\,+\\
+\,2\theta\,\frac{t^{2K}}{(2K)!}\,\frac{1}{M}\sum\limits_{N<n\le
N+M}\!\!e_{n}^{2K}\,=\\
=\,\sum\limits_{k =
0}^{K}\frac{(-1)^{k}t^{2k}}{(2k)!}\,\frac{(2k)!}{2^{k}k!}\biggl(1\,+\,\frac{0.2\theta_{1}(4A\sqrt[6\;]{3})^{k}}{\sqrt{L}}\biggr)\,+\,\theta_{2}
\sum\limits_{k =
1}^{K}\frac{|t|^{2k-1}}{(2k-1)!}\,\frac{(2k-1)!}{k!}\,\frac{1.6}{\sqrt{B}}\,\frac{(Be\pi^{2})^{k}}{\sqrt{L}}\,+\\
+\,2\theta_{3}\,\frac{t^{2K}}{(2K)!}\,\frac{(2K)!}{2^{K}K!}\biggl(1\,+\,\frac{0.2(4A\sqrt[6\;]{3})^{K}}{\sqrt{L}}\biggr),
\end{multline*}
where $\theta_{1} = 0$ for $k = 0$. After some transformations we
have
\[
f(t)\,=\,g(t)\,+\,\theta(r_{1}+r_{2}+r_{3}),
\]
where $g(t) = e^{-t^{2}/2}$,
\[
r_{1}\,=\,\frac{3}{K!}\biggl(\frac{t^{2}}{2}\biggr)^{K},\quad
r_{2}\,=\,\frac{0.4|t|}{\sqrt{L}}\sum\limits_{k =
1}^{K}\frac{(2A\sqrt[6\;]{3})^{k}}{k!}\,|t|^{2k-1},\quad
r_{3}\,=\,\frac{2.4|t|}{\sqrt{BL}}\sum\limits_{k =
1}^{K}\frac{(Be\pi^{2})^{k}}{k!}\,|t|^{2(k-1)}.
\]
Now let's consider the integral $I(\lambda)$,
\[
I(\lambda)\,=\,\int_{0}^{\lambda}\frac{|f(t)-g(t)|}{|t|}\,dt,
\]
where $\lambda > 1$. The we have the following estimate:
\begin{align*}
I(\lambda)\,& \le\,
\frac{3}{K}\,\frac{1}{K!}\biggl(\frac{\lambda^{2}}{2}\biggr)^{K}\,+\,\frac{0.4}{\sqrt{L}}\sum\limits_{k
=
1}^{K}\frac{(2A\sqrt[6\;]{3})^{k}}{k!}\,\frac{\lambda^{2k}}{2k}\,+\,\frac{2.4}{\sqrt{BL}}\sum\limits_{k
= 1}^{K}\frac{(Be\pi^{2})^{k}}{k!}\,\frac{\lambda^{2k-1}}{2k-1}\,<\\
&
<\,\biggl(\frac{\lambda^{2}e}{2K}\biggr)^{K}\,+\,\frac{0.2}{\sqrt{L}}\,\exp{\bigl(2A\sqrt[6\;]{3}\lambda^{2}\bigr)}\,+\,\frac{2.4}{\sqrt{BL}}\,\exp{\bigl(Be\pi^{2}\lambda^{2}\bigr)}\,<\\
&
<\,\biggl(\frac{\lambda^{2}e}{2K}\biggr)^{K}\,+\,\frac{0.5}{\sqrt{L}}\,\exp{\bigl(Be\pi^{2}\lambda^{2}\bigr)}.
\end{align*}
Setting
\[
\lambda\,=\,\frac{1}{2\pi}\sqrt{\frac{\ln{L}}{Be}},\quad
K\,=\,\biggl[\frac{Be}{2}\lambda^{2}\biggr]\,+\,1\,=\,\biggl[\frac{1}{8\pi^{2}}\,\ln{L}\biggr]\,+\,1,
\]
we get:
\begin{align*}
&
\frac{0.5}{\sqrt{L}}\,\exp{\bigl(Be\pi^{2}\lambda^{2}\bigr)}\,<\,\frac{0.5}{\sqrt[4\;]{L}},\\
&
\biggl(\frac{\lambda^{2}e}{2K}\biggr)^{K}\,\le\,B^{-K}\,\le\,\exp{\biggl(-\,\frac{\ln
B}{8\pi^{2}}\,\ln{L}\biggr)}\,<\,\exp{\biggl(-\,\frac{5}{\pi^{2}}\,\ln{L}\biggr)}\,<\,\frac{1}{\sqrt{L}}.
\end{align*}
Thus, $I(\lambda)<1/\sqrt[4\;]{L}$. Let $F(x)$ and $G(x)$ be the
distridution functions corresponding to the characteristic functions
$f(t)$ and $g(t)$, respectively. By the Berry\,-Esseen inequality
(see, for example, \cite[p. 20.3]{Loev_1962}), for any real $x$ we
have:
\[
|F(x)\,-\,G(x)|\,\le\,\frac{2}{\pi}I(\lambda)\,+\,\frac{12\sqrt{2}}{\pi\sqrt{\pi}}\,\lambda^{-1}\,<\,\frac{1}{\sqrt[4\;]{L}}\,+\,\frac{e^{20.5}\vep^{-1.5}}{\sqrt{\ln{L}}}
\,<\,\frac{A}{\sqrt{\ln{L}}}.
\]
Hence,
\[
F(x)\,=\,G(x)\,+\,\frac{\theta
A}{\sqrt{\ln{L}}}\,=\,\frac{1}{\sqrt{2\pi}}\int_{-\infty}^{\,x}\!e^{-u^{2}/2}du\,+\,\frac{\theta
A}{\sqrt{\ln{L}}}.
\]
The theorem is proved.

\vspace{0.2cm}

\textbf{ Corollary.} \emph{Let} $\Phi(x)$ \emph{be an arbitrary
positive function that increases unboundedly as} $x\to +\infty$.
\emph{Then the inequalities}
\[
\frac{1}{\Phi(n)}\,\frac{\sqrt{\ln\ln n}}{\ln
n}\,<\,|\gamma_{n}\,-\,t_{n}|\,\le\,\Phi(n)\,\frac{\sqrt{\ln\ln
n}}{\ln n}
\]
\emph{hold for almost all} $n$ \emph{in the following sense}:
\emph{the quantity of numbers} $n$, $n\le N$ \emph{that do not
satisfy these inequalities}, \emph{is} $o(N)$ \emph{as} $N\to
+\infty$.

\begin{center}
\textbf{\S 4. The Behavior of Quantities}
$\boldsymbol{S(\gamma_{n}+0)}$, $\boldsymbol{S(\gamma_{n}-0)}$
\textbf{``in the Mean''}
\end{center}

The function $S(t)$ is a piecewise smooth function with
discontinuities at the ordinates of zeros of $\zeta(s)$. It
decreases monotonically on every interval of discontinuity of the
form  $(\gamma_{n}, \gamma_{n+1})$. Thus, the right and left limits
of $S(t)$ at a points of discontinuity have an obvious geometric
sense: they are an upper and lower `peaks' of saw\,-tooth graph of
the function $S(t)$ (see fig. 1). The quantities $\Delta(n)$ are
connected with the values of $S(t)$ at a Gram's points ($\Delta(n) =
S(t_{n})$) and, similarly, the differences $\Delta_{n}$ are
connected with the values of $S(t)$ at the points of discontinuity
(lemma 7). This fact allows one to calculate the moments of the
quantities $S(\gamma_{n}\pm 0)$ and to find approximately the
distribution function for these quantities.

\begin{center}
\includegraphics{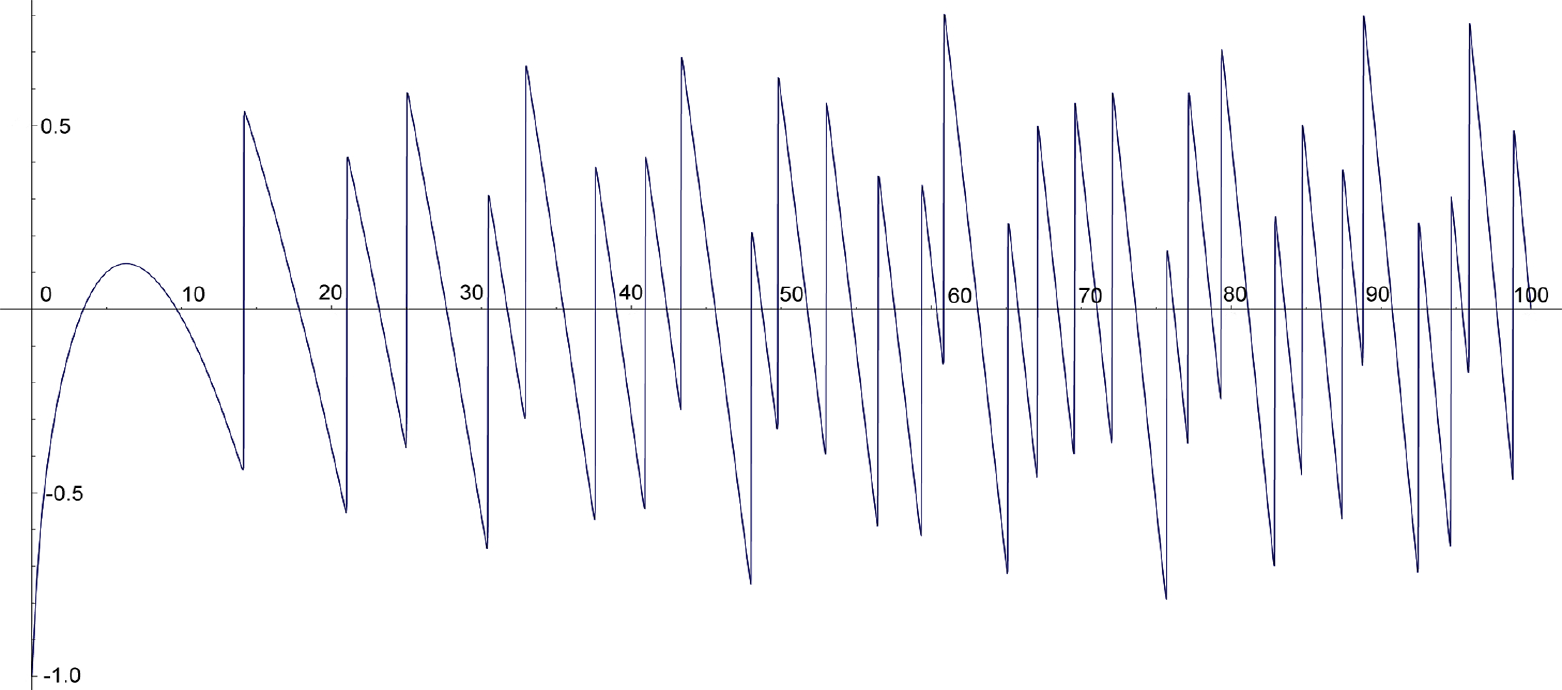}
\end{center}

\begin{center}
\fontsize{10}{10pt}\selectfont Fig. 1. The graph of the function
$S(t)$. Vertical segments on the graph correspond to the jumps of
$S(t)$ at a points of discontinuity, i.e. at ordinates of zeros of
$\zeta(s)$. \fontsize{12}{15pt}\selectfont
\end{center}

\vspace{0.2cm}

\textbf{Lemma 7.} \emph{Suppose} $n\to +\infty$. \emph{Then the
following equalities hold}
\[
S(\gamma_{n}+0)\,=\,-\Delta_{n}+\theta_{1}\kappa_{n}+O\biggl(\frac{\ln{n}}{n}\biggr),\quad
S(\gamma_{n}-0)\,=\,-\Delta_{n}-\theta_{2}\kappa_{n}+O\biggl(\frac{\ln{n}}{n}\biggr),
\]
\emph{where} $0\le \theta_{1},\theta_{2}\le 1$, \emph{and the
implied constants in} $O$'\emph{s are absolute}.

\vspace{0.2cm}

\emph{Proof.} By definition of $q_{n}$ we get $\gamma_{n} =
t_{n}+r_{n}$, where $r_{n} = q_{n}(t_{n+1}-t_{n})$. Applying
Riemann\,-von Mangoldt formula we obtain
\begin{multline*}
N(\gamma_{n}+0) =
\frac{1}{\pi}\vth(\gamma_{n})+1+S(\gamma_{n}+0)\,=\,\frac{1}{\pi}\vth(t_{n}+r_{n})+1+S(\gamma_{n}+0)\,=\\
=\,\frac{1}{\pi}\vth(t_{n})+1+\frac{r_{n}}{\pi}\vth'(t_{n})+\frac{r_{n}^{2}}{2\pi}\vth''(\xi)+1+S(\gamma_{n}+0),
\end{multline*}
where $\xi$ lies between $\gamma_{n}$ and $t_{n}$. By (\ref{Lab14}),
we have
\begin{align*}
&
\frac{r_{n}}{\pi}\vth'(t_{n})\,=\,q_{n}\,\frac{1}{\pi}\vth'(t_{n})(t_{n+1}-t_{n})\,=\,q_{n}\biggl(1\,+\,O\biggl(\frac{1}{n\ln
n}\biggr)\biggr)\,=\,q_{n}+O\biggl(\frac{1}{n}\biggr),\\
&
\frac{r_{n}^{2}}{2\pi}\vth''(\xi)\,=\,O\biggl(\frac{1}{\xi}\biggr)\,=\,O\biggl(\frac{\ln
n}{n}\biggr),
\end{align*}
and hence
\[
N(\gamma_{n}+0)\,=\,n\,+\,q_{n}\,+\,S(\gamma_{n}+0)\,+\,O\biggl(\frac{\ln
n}{n}\biggr).
\]
From the other hand, in the proof of theorem 1 we find that
\[
N(\gamma_{n}+0)\,=\,n\,+\,\theta_{n}(\kappa_{n}-1),\quad
0\le\theta_{n}\le 1.
\]
Comparing these expressions for $N(\gamma_{n}+0)$ and applying lemma
1 we get
\begin{align*}
&
S(\gamma_{n}+0)\,=\,-q_{n}+\theta_{n}(\kappa_{n}-1)+O\biggl(\frac{\ln
n}{n}\biggr)\,=\,-\Delta_{n}\,+\,\theta\,+\,\theta_{n}(\kappa_{n}-1)\,+\,O\biggl(\frac{\ln
n}{n}\biggr)\,=\\
& =\,-\Delta_{n}\,+\,\theta_{1}\kappa_{n}\,+\,O\biggl(\frac{\ln
n}{n}\biggr),\\
&
S(\gamma_{n}-0)\,=\,S(\gamma_{n}+0)\,-\,\kappa_{n}\,=\,-\Delta_{n}\,-\,\theta_{2}\kappa_{n}\,+\,O\biggl(\frac{\ln
n}{n}\biggr),
\end{align*}
where $0\le\theta_{1}\le 1$, $\theta_{2} = 1-\theta_{1}$. The lemma
is proved.

\vspace{0.2cm}

\textbf{Remark.} Lemma 7 allows one to estimate the constant $c$ in
the inequality of Theorem 1: $|\Delta_{n}|\le c\ln{n}$. Indeed,
using an upper bound for $|S(t)|$, from
\cite{Karatsuba_Korolev_2005} we obtain
\begin{align*}
|\Delta_{n}|\,=&
\,|S(\gamma_{n}+0)\,-\,\theta_{1}\kappa_{n}|\,+\,O\biggl(\frac{\ln
n}{n}\biggr)\,=\\
=&\,|(1-\theta_{1})S(\gamma_{n}+0)\,+\,\theta_{1}S(\gamma_{n}-0)|\,+\,O\biggl(\frac{\ln
n}{n}\biggr)\,\le\\
\le &
\,(1-\theta_{1})|S(\gamma_{n}+0)|\,+\,\theta_{1}|S(\gamma_{n}-0)|\,+\,O\biggl(\frac{\ln
n}{n}\biggr)\,\le\\
\le &
\,(1\,-\,\theta_{1}\,+\,\theta_{1})8.9\ln\gamma_{n}\,+\,O\biggl(\frac{\ln
n}{n}\biggr)<8.9\ln{n}.
\end{align*}
More general, suppose one has an estimation of the type $|S(t)|\le
cf(t)$ where the function $f(t)$ increases monotonically and $c$ is
a positive constant. Then by the same arguments we can conclude that
$|\Delta_{n}|\le(c+\vep)f(n)$ for any $\vep\!>\!0$ and $n\ge
n_{0}(\vep)$. Thus, if the Riemann hypothesis is true then the best
known bound for $|S(t)|$ (see \cite{Goldston_Gonek_2005}) implies
that
\[
|\Delta_{n}|\,\le\,\biggl(\frac{1}{2}+\vep\biggr)\frac{\ln n}{\ln\ln
n}.
\]

\vspace{0.2cm}

For the below, we need some assertions concerning the quantities
$\kappa_{n}$. The main purpose of these assertions is to show that
these quantities behave as a constants `in the mean'.

\vspace{0.2cm}

\textbf{Definition 8. }Let $j\ge 1$ be an integer. Denote by
$n_{j}(T)$ the number of ordinates $\gamma$ of zeros of $\zeta(s)$,
$0<\gamma\le T$, whose multiplicities are equal to $j$.

\vspace{0.2cm}

First we show that the number $n_{j}(T)$ is sufficiently small for
large $j$. In essence, the proof of this fact repeats
word\,-for\,-word an upper estimate for the number of zeros of
$\zeta(s)$ with given multiplicity. This proof is based on two
facts: 1) an absolute value of the difference $S(t+h)-S(t)$ is very
large whet $t$ is close to the ordinate $\gamma$ with high
multiplicity for $h\asymp (\ln t)^{-1}$ (lemma 8); 2) the measure of
the subset of given interval where this difference is large, is
small enough (lemma 9).

In what follows, $T\ge T_{0}(\vep)>0$, $h =
\frac{\pi}{3}\bigl(\ln{\frac{T}{2\pi}}\bigr)^{-1}$. The parameter
$H$ is supposed to be the following. The below lemma is proved in
\cite{Korolev_2006} only for $H = T^{27/82+\vep}$. This restriction
is inconvenient in some cases. The proof of lemma 9 is based on the
density theorem for the zeros of the Riemann zeta\,-function, that
is, on the upper bound for the number
\begin{equation}
N(\sigma,T+H)\,-\,N(\sigma,T)\label{Lab29}
\end{equation}
of zeros of $\zeta(s)$ in the region $\sigma < \RRe s <1$, $T<\IIm
s\le T+H$ (see \cite{Karatsuba_Korolev_2006}). However, the proof of
this bound is valid also in the case $T^{27/82+\vep_{1}} \le H \le
T^{27/82+\vep_{2}}$, where $\vep_{1}$, $\vep_{2}$ are an arbitrary
numbers with the condition $0.9\vep \le \vep_{1}<\vep_{2}\le \vep$.
Therefore we suppose that $H$ satisfies the above conditions with
$\vep_{1} = 0.9\vep$, $\vep_{2} = \vep$.

\vspace{0.2cm}

\textbf{Lemma 8.} \emph{Let} $\gamma$ \emph{be the imaginary part of
a zero of} $\zeta(s)$ \emph{such that} $T\le \gamma - h <\gamma \le
T+H$. \emph{Suppose that the interval} $(\gamma - h,\gamma]$
\emph{contains exactly} $m$ \emph{ordinates of zeros of} $\zeta(s)$
\emph{whose multiplicities are equal to} $j$. \emph{Then we have}
\[
S(t+2h)\,-\,S(t)\,\ge\,mj\,-\,0.5
\]
\emph{for any} $t$ \emph{in} $(\gamma-2h, \gamma-h]$.

\emph{Proof.} Suppose that $\gamma - 2h < t<\gamma - h$ and that
$\gamma^{(1)}<\ldots <\gamma^{(r)}$ are all distinct ordinates of
zeros of $\zeta(s)$ contained in $(t,t+2h]$. Using the same
arguments as above (see theorem 2) we obtain the identity
\begin{multline}
S(t+2h)-S(t)\,=\,\bigl(S(t+2h)-S(\gamma^{(r)}+0)\bigl)\,+\,\sum\limits_{\nu
= 1}^{r}\bigl(S(\gamma^{(\nu)}+0)\,-\,S(\gamma^{(\nu)}-0)\bigr)\,+\\
+\,\sum\limits_{\nu =
2}^{r}\bigl(S(\gamma^{(\nu)}-0)\,-\,S(\gamma^{(\nu-1)}+0)\bigr)\,+\,\bigl(S(\gamma^{(1)}-0)\,-\,S(t)\bigr).\label{Lab30}
\end{multline}
All the differences $S(\gamma^{(\nu)}+0)\,-\,S(\gamma^{(\nu)}-0) =
\kappa(\gamma^{(\nu)})$ are positive. Moreover, the inequalities
$\gamma - h > t$, $\gamma < t+2h$ imply that the interval $(\gamma -
h,\gamma]$ is contained entirely in $(t,t+2h]$, and there are at
least $m$ ordinates among $\gamma^{(1)},\ldots, \gamma^{(r)}$, whose
multiplicities are equal to $j$. Hence,
\[
\sum\limits_{\nu =
1}^{r}\bigl(S(\gamma^{(\nu)}+0)\,-\,S(\gamma^{(\nu)}-0)\bigr)\,\ge\,mj.
\]
Transforming all other differences in the right\,-hand side of
(\ref{Lab30}) by Lagrange's mean value theorem and using the
relation
\[
S'(t)\,=\,-\frac{1}{2\pi}\ln\frac{t}{2\pi}\,+\,O(t^{-2})\,=\,-\biggl(\frac{1}{6h}\,+\,\frac{\theta}{\pi}\frac{H}{T}\biggr),
\]
valid for $T\le t\le T+H$, we obtain
\begin{multline*}
S(t+2h)\,-\,S(t)\,\ge\,mj\,-\,\bigl((t+2h-\gamma^{(r)})\,+\,\sum\limits_{\nu
= 2}^{r}(\gamma^{(\nu)}\,-\,\gamma^{(\nu -
1)})\,+\,(\gamma^{(1)}-t)\bigr)\biggl(\frac{1}{6h}\,+\,\frac{H}{\pi
T}\biggr)\,=\\
=\,mj\,-\,2h\biggl(\frac{1}{6h}\,+\,\frac{H}{\pi
T}\biggr)\,>\,mj\,-\,0.5.
\end{multline*}
The lemma is proved.

\vspace{0.2cm}

\textbf{Lemma 9. }\emph{Let} $D(\lambda)$ \emph{be the set of
points} $t$ \emph{such that} $T\le t\le T+H$ \emph{and}
$|S(t+2h)\,-\,S(t)|\,\ge\,\lambda$. \emph{Then the inequality}
\[
\text{mes}\,D(\lambda)\,\le\,e^{4}He^{-C\lambda},
\]
\emph{holds for any} $\lambda$ \emph{with} $C = \pi
e\sqrt{\frac{2}{5}e}\,A^{-1}$.

\vspace{0.2cm}

For the proof, see \cite{Korolev_2006}.

\vspace{0.2cm}

\textbf{Lemma 10.} \emph{For any} $j\ge 1$, \emph{the following
estimation holds}:
\[
n_{j}(T+H)\,-\,n_{j}(T)\,<\,e^{7.2}\bigl(N(T+H)\,-\,N(T)\bigr)e^{-Cj},
\]
\emph{where the constant} $C$ \emph{is defined in lemma 9}.

\vspace{0.2cm}

\emph{Proof.} We follow the proof of theorems B and C from
\cite{Korolev_2006}. Let $j$ be the multiplicity of an ordinate
$\gamma$, $T<\gamma \le T+H$. Then the inequality $|S(t)|<8.9\ln{t}$
implies
\[
j\,=\,S(\gamma +
0)\,-\,S(\gamma-0)\,\le\,17.8\ln\gamma\,<\,18\ln{T}.
\]
Hence, the difference $n_{j}(T+H)-n_{j}(T)$ equals to zero for $j\ge
18\ln{T}$. In the case $1\le j\le C^{-1}\ln{2}$ we obviously have
$1\le 2e^{-Cj}$ and therefore
\[
n_{j}(T+H)\,-\,n_{j}(T)\le N(T+H)-N(T)\le 2(N(T+H)-N(T))e^{-Cj}.
\]
Thus, it is sufficient to consider the case
\[
C^{-1}\ln{2}\,<\,j\,<\,18\ln{T}.
\]
Suppose $\gamma_{(1)}$ is the largest ordinate with multiplicity
equal to $j$ in the interval $(T,T+H]$. Then we denote $E_{1} =
(\gamma_{(1)}-h,\gamma_{(1)}]$ and stand the symbol $\gamma_{(2)}$
for the largest ordinate whose multiplicity equals to $j$ in the
interval $(T,\gamma_{(1)}-h]$. Further, we denote $E_{2} =
(\gamma_{(2)}-h,\gamma_{(2)}]$ and stand the symbol $\gamma_{(3)}$
for the largest ordinate whose multiplicity equals to $j$ in the
interval $(T,\gamma_{(2)}-h]$ and so on.

We continue this construction until there are no such ordinates in
the interval $(T,\gamma_{(r)}-h]$ or until we find such ordinate
$\gamma_{(r)}$ satisfying the condition $\gamma_{(r)}-h<T\le
\gamma_{(r)}$. In both cases we set $E_{r} =
(\gamma_{(r)}-h,\gamma_{(r)}]$.

The intervals $E_{1}, E_{2}, \ldots, E_{r}$  are pairwise disjoint
and have the same length $h$. Moreover, their union contains all
ordinates of zeros with multiplicity equals to $j$ lying in
$(T,T+H]$. Now we partition the intervals constructed into classes
$\mathcal{E}_{1}, \mathcal{E}_{2}, \ldots$ by putting into class
$\mathcal{E}_{m}$ the intervals containing exactly $m$ of desired
ordinates. If $k_{m}$ is the number of intervals in the class
$\mathcal{E}_{m}$ then
\begin{equation}
n_{j}(T+H)\,-\,n_{j}(T)\,\le\,1\cdot k_{1}\,+\,2\cdot k_{2}\,+\ldots
+ mk_{m}\,+\ldots \label{Lab31}
\end{equation}
Let us find an upper bound for each of quantities $k_{m}$. Suppose
the interval $E = (\gamma - h,\gamma]$ belongs to $\mathcal{E}_{m}$.
We then set $E' = E-h = (\gamma-2h, \gamma - h]$. By lemma 8,
\[
S(t+2h)\,-\,S(t)\,\ge\,mj\,-\,0.5
\]
for any $t\in E'$. Since $E'$ is contained in $(T-2h,T+H]$ then $E'$
is entirely contained in $D(mj-0.5)$ where $D(\lambda)$ denotes the
set of points $t\in (T-2h,T+H]$ satisfying the inequality
$|S(t+2h)-S(t)|\ge\lambda$. After carrying out this construction for
any interval $E$ in $\mathcal{E}_{m}$, we see that all the $k_{m}$
intervals are contained in $D(mj-0.5)$. Since they are pairwise
disjoint and have the same length $h$, their total length does not
exceed the measure of $D(mj-0.5)$, that is
\[
k_{m}h\,\le\,\text{mes}\,D(mj-0.5).
\]
By lemma 9, we get
\begin{align*}
&
\text{mes}\,D(mj-0.5)\,\le\,e^{4}(H+2h)e^{-C(mj-0.5)}\,<\,e^{4+C}He^{-Cmj},\\
& k_{m}\,\le\,e^{4+C}Hh^{-1}e^{-Cmj}.
\end{align*}
Returning to (\ref{Lab31}), we obtain
\[
n_{j}(T+H)\,-\,n_{j}(T)\,\le\,e^{4+C}Hh^{-1}\sum\limits_{m =
1}^{+\infty}me^{-Cmj}\,=\,e^{4+C}Hh^{-1}\,\frac{e^{-Cj}}{(1-e^{-Cj})^{2}}.
\]
Since $e^{-Cj}\le \tfrac{1}{2}$, then
\[
n_{j}(T+H)\,-\,n_{j}(T)\,\le\,4e^{4+C}Hh^{-1}e^{-Cj}\,<\,e^{7.2}\bigl(N(T+H)-N(T)\bigr)e^{-Cj}.
\]
The lemma is proved.

\vspace{0.2cm}

\textbf{Corollary.} \hfill \emph{Let} \hfill $\nu(j) =
\nu(j;M,N)$\hfill \emph{be \hfill the \hfill number \hfill of \hfill
distinct \hfill ordinates} \hfill $\gamma_{n}$, \\ $N<n\le N+M$,
\emph{whose multiplicities are equal to} $j$. \emph{Then the
following estimation holds}: $\nu(j)\le e^{7.3}Me^{-Cj}$.

\vspace{0.2cm}

\emph{Proof.} Setting $T = \gamma_{N}$, $H =
\gamma_{N+M}-\gamma_{N}$ in lemma and noting that the multiplicities
of $\gamma_{N}$, $\gamma_{N+M}$ are of order $\ln{N}$, we obtain
\[
\nu(j)\,\le\,e^{7.2}(M+O(\ln N))e^{-Cj}\,<\,e^{7.3}Me^{-Cj}.
\]

The bounds of lemma 10 and it's corollary allows one to study the
behavior of multiplicities $\kappa_{n}$ `in the mean'.

\vspace{0.2cm}

\textbf{Definition 9.} Suppose $a>0$. We denote by
\[
K_{0}(a)\,=\,K_{0}(a;M,N)\,=\,\sum\limits_{N<n\le N+M}\kappa_{n}^{a}
\]
the sum over \emph{all} ordinates $\gamma_{n}$, $N < n\le N+M$.
Further, we denote by
\[
K_{1}(a)\,=\,K_{1}(a;M,N)\,=\,\prsum\limits_{N<n\le
N+M}\kappa_{n}^{a}
\]
the sum over \emph{all distinct} ordinates $\gamma_{n}$, $N<n\le
N+M$.

Thus, if (\ref{Lab9}) holds then the sum $K_{0}(a)$ contains, for
example, all the terms $\kappa_{l}^{a},$
$\kappa_{l+1}^{a},\ldots,\kappa_{l+p-1}^{a}$, though the sum
$K_{1}(a)$ contains only the term $\kappa_{l}^{a}$. It's easy to see
that
\[
K_{0}(a)\,=\,\prsum\limits_{N<n\le N+M}\kappa_{n}\cdot
\kappa_{n}^{a}\,=\,K_{1}(a+1).
\]

\textbf{Lemma 11.} \emph{Suppose} $a\ge 1$. \emph{The the following
estimation holds}
\[
K_{1}(a)\,<\,e^{8.4}\,\frac{\Gamma(a+1)}{C^{a+1}}M,
\]
\emph{where the constant} $C$ \emph{is defined in lemma 9}.

\vspace{0.2cm}

\emph{Proof.} Using the assertion and the notations of lemma 10 and
it's corollary, we find
\begin{equation}
K_{1}(a) = \sum\limits_{j\ge
1}j^{a}\nu(j)\,\le\,e^{7.3}M\sum\limits_{j =
1}^{\infty}j^{a}e^{-Cj}.\label{Lab32}
\end{equation}
The function $y(x) = x^{a}e^{-Cx}$ increases on the segment $1\le
x\le x_{a}$, $x_{a} = aC^{-1}$, up to it's maximum
\begin{equation}
\biggl(\frac{a}{C}\biggr)^{a}e^{-a}\,=\,\biggl(\frac{a}{Ce}\biggr)^{a}\label{Lab33}
\end{equation}
at a point $x_{a}$ and then decreases monotonically. Let's define
the integer $m$ from the inequalities $m<x_{a}\le m+1$. Then we
estimate the terms of the sum (\ref{Lab32}) corresponding to $j =
m,m+1$ by the quantity (\ref{Lab33}) and all other terms by the
integrals of the function $y(x)$. Thus we get
\begin{multline}
\sum\limits_{j =
1}^{\infty}j^{a}e^{-Cj}\,\le\,\biggl(\int_{0}^{m}\,+\,\int_{m+1}^{\infty}\biggr)x^{a}e^{-Cx}dx\,+\,2\biggl(\frac{a}{Ce}\biggr)^{a}\,<\\
<\,\int_{0}^{\infty}x^{a}e^{-Cx}dx\,+\,2\biggl(\frac{a}{Ce}\biggr)^{a}\,<\,\frac{\Gamma(a+1)}{C^{a+1}}\,+\,2\biggl(\frac{a}{Ce}\biggr)^{a}\,<\,\frac{3\Gamma(a+1)}{C^{a+1}}.
\end{multline}
Therefore, the desired assertion follows from the last inequality.

\vspace{0.2cm}

The estimate of lemma 11 together with the equalities of lemma 7
allow one to calculate the moments of the quantities
$S(\gamma_{n}\pm 0)$.

\vspace{0.2cm}

\textbf{Theorem 11.} \emph{Let} $k$ \emph{be an integer with the
condition} $1\le k\le \sqrt{L}$. \emph{Then the following relation
holds}
\begin{align*}
& \sum\limits_{N<n\le
N+M}S^{2k}(\gamma_{n}+0)\,=\,\frac{\vk(2k)}{(2\pi)^{2k}}\,ML^{k}\bigl(1\,+\,19\theta_{1} A(4A\sqrt{3})^{k}L^{-0.5}\bigr),\\
& \sum\limits_{N<n\le
N+M}S^{2k}(\gamma_{n}-0)\,=\,\frac{\vk(2k)}{(2\pi)^{2k}}\,ML^{k}\bigl(1\,+\,19\theta_{2}
A(4A\sqrt{3})^{k}L^{-0.5}\bigr).
\end{align*}

\vspace{0.2cm}

\emph{Proof.} Using the same arguments as above (see the proof ot
theorem 7) and the equality of lemma 7 we get
\begin{align*}
& S(\gamma_{n}+0)\,=\,-\,\Delta_{n}\,+\,1.1\theta\kappa_{n},\\
& S^{2k}(\gamma_{n}+0)\,=\,\Delta_{n}^{2k}\,+\,\theta\,k
2^{2k}\bigl(|\Delta_{n}|^{2k-1}\kappa_{n}\,+\,(1.1)^{2k}\kappa_{n}^{2k}\bigr),
\end{align*}
and finally
\[
\sum\limits_{N<n\le
N+M}S^{2k}(\gamma_{n}+0)\,=\,\Sigma_{1}\,+\,\theta\,k
2^{2k}\bigl((1.1)^{2k}\Sigma_{2}\,+\,\Sigma_{3}\bigr),
\]
where
\begin{align*}
& \Sigma_{1}\,=\,\sum\limits_{N<n\le N+M}\Delta_{n}^{2k},\quad
\Sigma_{2}\,=\,\sum\limits_{N<n\le
N+M}\kappa_{n}^{2k}\,=\,\prsum\limits_{N<n\le
N+M}\kappa_{n}^{2k+1},\\
& \Sigma_{3}\,=\,\sum\limits_{N<n\le
N+M}|\Delta_{n}|^{2k-1}\kappa_{n}.
\end{align*}
First, by lemma 11 we conclude that
\[
\Sigma_{2}\,\le\,\frac{e^{8.4}}{C^{2k+2}}\,(2k+1)!
M\,=\,\frac{\vk(2k)}{(2\pi)^{2k}}ML^{k}\delta_{2},
\]
where
\begin{multline*}
\delta_{2} =
\frac{e^{8.4}}{C^{2}}\,\biggl(\frac{2\pi}{C}\biggr)^{2k}\frac{k!(2k+1)!}{(2k)!L^{k}}\,\le\,\frac{e^{8.4}}{C^{2}}\biggl(\frac{4\pi^{2}}{C^{2}L}\biggr)^{k}
(2k+1)\,2\sqrt{k}\biggl(\frac{k}{e}\biggr)^{k}\,<\\
<\,\frac{e^{11}}{C^{2}}\,k\sqrt{k}\biggl(\frac{4\pi^{2}k}{C^{2}L
e}\biggr)^{k}\,\le\,\frac{e^{11}}{C^{2}}\biggl(\frac{4\pi^{2}\sqrt{3}k}{C^{2}L
e}\biggr)^{k}.
\end{multline*}
The right\,-hand side has its maximum at the point $k = 1$ for $1\le
k\le \sqrt{L}$. Hence,
\[
\delta_{2}\,\le\,\frac{e^{11}}{C^{2}}\,\frac{4\pi^{2}\sqrt{3}}{C^{2}L}\,<\,\frac{1}{\sqrt{L}}.
\]
Moreover, we have
\[
\delta_{2}^{1/(2k)}\,=\,\biggl(\frac{e^{11}}{C^{2}}\biggr)^{1/(2k)}\sqrt{\frac{4\pi^{2}\sqrt{3}k}{C^{2}Le}}\,\le\,
\sqrt{\frac{e^{11}}{C^{2}}}\,\frac{2\pi\sqrt[4\;]{3}}{C\sqrt{e}}\,\sqrt{\frac{k}{L}}\,<\,e^{2.8}A^{2}\sqrt{\frac{k}{L}}.
\]
Further, by the corollary of theorem 4 we obtain
\[
\Sigma_{1}\,\le\,\frac{\vk(2k)}{(2\pi)^{2k}}\,ML^{k}\delta_{1},\quad
\delta_{1}\,=\,1.1A^{k-1}.
\]
Thus
\begin{align*}
& k\,2^{2k}\Sigma_{3}\,\le\,k\,2^{2k}\,\Sigma_{1}^{1-1/(2k)}\Sigma_{2}^{1/(2k)}\,\le\,\frac{\vk(2k)}{(2\pi)^{2k}}\,ML^{k}\delta_{3},\\
&
\delta_{3}\,=\,k\,2^{2k}\delta_{1}^{1-1/(2k)}\delta_{2}^{1/(2k)}\,\le
k2^{2k}\biggl(\frac{1.1}{A}\biggr)^{1-1/(2k)}A^{k-1/2}\,e^{2.8}A^{2}\sqrt{\frac{k}{L}}\,\le\,\\
&\le\,\sqrt{1.1}e^{2.8}A(4A)^{k}\frac{k\sqrt{k}}{\sqrt{L}}\,<\,\frac{18A(4\sqrt{3}A)^{k}}{\sqrt{L}}.
\end{align*}
Finally we get
\[
\sum\limits_{N<n\le
N+M}S^{2k}(\gamma_{n}+0)\,=\,\frac{\vk(2k)}{(2\pi)^{2k}}\,ML^{k}(1\,+\,\theta\delta),
\]
where
\[
\delta\,=\,\frac{1.1A^{k}}{\sqrt{L}}\,+\,\frac{18A(4A\sqrt{3})^{k}}{\sqrt{L}}\,+\,\frac{k\,2^{2k}}{\sqrt{L}}\,<\,\frac{19A(4A\sqrt{3})^{k}}{\sqrt{L}}.
\]
The second equality can be proved in a similar way. The theorem is
proved.

\vspace{0.2cm}

\textbf{Theorem 12.} \emph{Let} $k$ \emph{be an integer with the
condition} $1\le k\le \sqrt{L}$. \emph{Then the following relation
holds}
\[
\biggl|\sum\limits_{N<n\le
N+M}S^{2k-1}(\gamma_{n}+0)\biggr|<c_{k}ML^{k-1},\quad
\biggl|\sum\limits_{N<n\le
N+M}S^{2k-1}(\gamma_{n}-0)\biggr|<c_{k}ML^{k-1},
\]
\emph{where the quantity} $c_{k}$ \emph{can be set to be equal to}
$(eA)^{3}$ \emph{for} $k = 1$ \emph{and to} $3.7B^{-0.5}(Bk)^{k}$,
\emph{for} $k\ge 2$.

\vspace{0.2cm}

\emph{Proof.} It's sufficient to estimate the first sum. Suppose $k
= 1$. By lemmas 7,11 and theorem 5 we have:
\begin{multline*}
\biggl|\sum\limits_{N<n\le
N+M}S(\gamma_{n}+0)\biggr|\,\le\,\biggl|\sum\limits_{N<n\le
N+M}\Delta_{n}\biggr|\,+\,1.01\sum\limits_{N<n\le
N+M}\!\!\kappa_{n}\,\le\,3.6\sqrt{B}M\,+\,e^{2.7}A^{3}M\,<\\
<\,(e A)^{3}M.
\end{multline*}
Suppose now $k\ge 2$. Summing both parts of the equality
\[
S^{2k-1}(\gamma_{n}+0)\,=\,-\Delta_{n}^{2k-1}\,+\,0.5\theta
k\,2^{2k}\bigl((1.01)^{2k-1}\kappa_{n}^{2k-1}\,+\,\Delta_{n}^{2k-2}\kappa_{n}\bigr),
\]
we get
\[
\biggl|\sum\limits_{N<n\le N+M}S^{2k-1}(\gamma_{n}+0)\biggr|\,\le\,
\biggl|\sum\limits_{N<n\le
N+M}\Delta_{n}^{2k-1}\biggr|\,+\,0.5k\,2^{2k}\bigl((1.01)^{2k-1}\Sigma_{1}\,+\,\Sigma_{2}\bigr),
\]
where
\[
\Sigma_{1}\,=\,\sum\limits_{N<n\le
N+M}\kappa_{n}^{2k-1}\,=\,\prsum\limits_{N<n\le
N+M}\kappa_{n}^{2k},\quad \Sigma_{2}\,=\,\sum\limits_{N<n\le
N+M}\Delta_{n}^{2k-2}\kappa_{n}.
\]
By lemma 11,
\begin{align*}
&
\Sigma_{1}\,\le\,\frac{e^{8.4}}{C^{2k+1}}\,(2k)!M\,<\,e^{7.6}A\sqrt{k}\biggl(\frac{2k}{Ce}\biggr)^{2k}M,\\
&
0.5k\,2^{2k}(1.01)^{2k-1}\Sigma_{1}\,<\,e^{7}A\biggl(\frac{2k}{C}\biggr)^{2k}M\,=\,k^{k}ML^{k-1}\delta_{1},\\
& \delta_{1}\,=\,e^{7}AL
\biggl(\frac{4k}{C^{2}L}\biggr)^{k}\,\le\,e^{7}AL\biggl(\frac{8}{C^{2}L}\biggr)^{2}\,<\,\frac{1}{\sqrt{L}}\,<\,\frac{0.05}{\sqrt{B}}.
\end{align*}
Further, applying H\"{o}lder's inequality to the sum $\Sigma_{2}$ we
obtain $\Sigma_{2}\,\le\,\Sigma_{3}^{1-1/k}\Sigma_{4}^{1/k}$, where
\[
\Sigma_{3}\,=\,\sum\limits_{N<n\le N+M}\Delta_{n}^{2k},\quad
\Sigma_{4}\,=\,\sum\limits_{N<n\le
N+M}\kappa_{n}^{k}\,=\,\prsum\limits_{N<n\le N+M}\kappa_{n}^{k+1}.
\]
The corollary of theorem 4 and the estimate of lemma 1 imply
together
\begin{align*}
&
\Sigma_{3}\,\le\,\frac{1.1}{A}\,\frac{\vk(2k)}{(2\pi)^{2k}}\,M(AL)^{k}\,<\,\frac{\sqrt{e}}{A}\biggl(\frac{kAL}{\pi^{2}e}\biggr)^{k}M,\\
&
\Sigma_{4}\,\le\,\frac{e^{8.4}}{C^{k+2}}\,(k+1)!\,M\,<\,e^{6}A^{2}k\sqrt{k}\biggl(\frac{k}{Ce}\biggr)^{k}\!\!M.
\end{align*}
Hence
\begin{multline*}
0.5k\,2^{2k}\Sigma_{2}\,\le\,0.5k\,2^{2k}\,\frac{\sqrt[4\;]{e}}{\sqrt{A}}\biggl(\frac{kAL}{\pi^{2}e}\biggr)^{k-1}\!\!
e^{3}A\sqrt{3}\,\frac{k}{Ce}M\,<\\
<\,25\sqrt{A}\bigl(\tfrac{1}{4}kA\bigr)^{k}ML^{k-1}\,<\,\frac{0.05}{\sqrt{B}}\,(Bk)^{k}ML^{k-1}.
\end{multline*}
Now the desired assertion follows from theorem 5:
\[
\biggl|\sum\limits_{N<n\le
N+M}\!S^{2k-1}(\gamma_{n}+0)\,\biggr|\,<\,\frac{(Bk)^{k}}{\sqrt{B}}\,ML^{k-1}(3.6\,+\,0.05\,+\,0.05)\,=\,\frac{3.7}{\sqrt{B}}\,(Bk)^{k}ML^{k-1}.
\]
The theorem is proved.

\vspace{0.2cm}

\textbf{Theorem 13.} \emph{Let} $a$ \emph{be an arbitrary real
number with the condition}
\[
0<a<\frac{e\ln{L}}{10A\ln\ln L}.
\]
\emph{Then the following asymptotic formulae hold}
\begin{align*}
& \sum\limits_{N<n\le
N+M}|S(\gamma_{n}+0)|^{a}\,=\,\frac{\vk(a)}{(2\pi)^{a}}\,ML^{0.5a}\bigl(1\,+\,\theta_{1}\delta\bigr),\\
& \sum\limits_{N<n\le
N+M}|S(\gamma_{n}-0)|^{a}\,=\,\frac{\vk(a)}{(2\pi)^{a}}\,ML^{0.5a}\bigl(1\,+\,\theta_{2}\delta\bigr),
\end{align*}
\emph{where} $\delta = c(A^{-1}\ln{L})^{-\gamma}$. \emph{Therefore,
the constants} $c$ \emph{and} $\gamma$ \emph{can be set to be equal
to} $e^{5}A^{3}$ \emph{and} $0.5a$ \emph{if} $0<a\le 1$, \emph{and
to} $2^{15.5}$ \emph{and} $\tfrac{1}{2}+\bigl\{\frac{a-1}{2}\bigr\}$
\emph{if} $a>1$. \emph{In particular}, \emph{in the case} $a = 1$
\emph{we get}
\begin{align*}
& \sum\limits_{N<n\le
N+M}|S(\gamma_{n}+0)|\,=\,\frac{M}{\pi\sqrt{\pi}}\,\sqrt{\ln\ln N}\bigl(1\,+\,O\bigl((\ln\ln\ln N)^{-0.5}\bigr)\bigr),\\
& \sum\limits_{N<n\le
N+M}|S(\gamma_{n}-0)|\,=\,\frac{M}{\pi\sqrt{\pi}}\,\sqrt{\ln\ln
N}\bigl(1\,+\,O\bigl((\ln\ln\ln N)^{-0.5}\bigr)\bigr).
\end{align*}

\vspace{0.2cm}

\emph{Proof.} Suppose first that $0<a\le 1$. Using lemma 7 and
considering the cases $|\Delta_{n}|\le 2\kappa_{n}$ and
$|\Delta_{n}|>2\kappa_{n}$ separately we obtain
\[
|S(\gamma_{n}+0)|^{a}\,=\,|\Delta_{n}|^{a}\,+\,3.01\theta\kappa_{n}.
\]
Summing both parts over $n$ we get
\[
\sum\limits_{N<n\le
N+M}|S(\gamma_{n}+0)|^{a}\,=\,\Sigma_{1}\,+\,3.01\theta\Sigma_{2},
\]
where
\[
\Sigma_{1}\,=\,\sum\limits_{N<n\le N+M}|\Delta_{n}|^{a},\quad
\Sigma_{2}\,=\,\sum\limits_{N<n\le
N+M}\kappa_{n}\,=\,\prsum\limits_{N<n\le N+M}\kappa_{n}^{2}.
\]
Applying lemma 11 to the sum $\Sigma_{2}$ we have
\begin{align*}
&
\Sigma_{2}\,\le\,\frac{2e^{8.4}}{C^{3}}\,M\,=\,\frac{\vk(a)}{(2\pi)^{a}}\,ML^{0.5a}\delta_{2},\\
&
\delta_{2}\,\le\,\frac{\pi^{a+1/2}}{\Gamma\bigl(\frac{a+1}{2}\bigr)}\,\frac{2e^{8.4}}{C^{3}L^{0.5a}}\,\le\,2\pi^{1.5}e^{8.4}C^{-3}L^{-0.5a}\,<\,0.5e^{5}A^{3}L^{-0.5}.
\end{align*}
The application of asymptotic formula for the sum $\Sigma_{1}$
yields
\[
\sum\limits_{N<n\le
N+M}|S(\gamma_{n}+0)|^{a}\,=\,\frac{\vk(a)}{(2\pi)^{a}}\,ML^{0.5a}\bigl(1\,+\,\theta\delta\bigr),
\]
where
\[
\delta\,\le\,91(A^{-1}\ln{L})^{-0.5a}\,+\,0.5e^{5}A^{3}L^{-0.5a}\,<\,e^{5}A^{3}(A^{-1}\ln{L})^{-0.5a}.
\]
Suppose now $a>1$. By easy\,-to\,-check relation
\[
|S(\gamma_{n}+0)|^{a}\,=\,|\Delta_{n}|^{a}\,+\,(3.01)^{a}\theta\bigl(\kappa_{n}^{a}\,+\,|\Delta_{n}|^{a-1}\kappa_{n}\bigr)
\]
we obtain the following expression for the initial sum:
\[
\sum\limits_{N<n\le
N+M}|S(\gamma_{n}+0)|^{a}\,=\,\Sigma_{1}\,+\,(3.01)^{a}\theta\bigl(\Sigma_{2}\,+\,\Sigma_{1}^{1-1/a}\Sigma_{2}^{1/a}\bigr),
\]
where
\[
\Sigma_{2}\,=\,\sum\limits_{N<n\le
N+M}\kappa_{n}^{a}\,=\,\prsum\limits_{N<n\le N+M}\kappa_{n}^{a+1},
\]
(the symbol $\Sigma_{1}$ denotes the same sum as above). Using again
the estimation of lemma 11 together with the duplication formula for
the gamma\,-function we obtain
\begin{align*}
&
\Sigma_{2}\,\le\,\frac{e^{8.4}}{C^{a+2}}\,\Gamma(a+2)M\,\le\,\frac{\vk(a)}{(2\pi)^{a}}\,ML^{0.5a}\delta_{1},\\
&
\delta_{1}\,=\,\frac{\pi^{a+1/2}}{\Gamma\bigl(\frac{a+1}{2}\bigr)}\,\frac{\Gamma(a+2)}{C^{a+2}}\,\frac{e^{8.4}}{L^{0.5a}}\,<\,e^{9}C^{-2}\,
\frac{\Gamma(a+2)}{\Gamma\bigl(\frac{a+1}{2}\bigr)}\,\biggl(\frac{\pi}{C\sqrt{L}}\biggr)^{a}\,=
\end{align*}
\begin{align}
&
=\,e^{9}C^{-2}(a+2)\Gamma\biggl(\frac{a}{2}+1\biggr)\biggl(\frac{2\pi}{C\sqrt{L}}\biggr)^{a}\,\le\,e^{9}C^{-2}\biggl(\frac{4\pi
a}{C\sqrt{L}}\biggr)^{a}.\label{Lab34}
\end{align}
Since $a = o(\ln L)$, we obviously have
\[
\delta_{1}\,\le\,e^{9}C^{-2}\biggl(\frac{1}{\sqrt[3\;]{L}}\biggr)^{a}\,<\,\frac{1}{\sqrt[4\;]{L}}.
\]
Next, noting that
\[
\Sigma_{1}\,<\,1.01\,\frac{\vk(a)}{(2\pi)^{a}}\,ML^{0.5a}
\]
and applying the bound (\ref{Lab34}) for $\delta_{1}$, we obtain
\[
\Sigma_{1}^{1-1/a}\Sigma_{2}^{1/a}\,<\,\frac{\vk(a)}{(2\pi)^{a}}\,ML^{0.5a}\,\delta_{2},\quad
\delta_{2}\,<\,1.01e^{9}C^{-2}\,\frac{4\pi
a}{C\sqrt{L}}\,<\,\frac{1}{\sqrt[4\;]{L}}.
\]
Thus,
\[
\sum\limits_{N<n\le
N+M}|S(\gamma_{n}+0)|^{a}\,=\,\frac{\vk(a)}{(2\pi)^{a}}\,ML^{0.5a}\bigl(1\,+\,\theta\delta\bigr),
\]
where
\[
\delta =
2^{15.5}\bigl(A^{-1}\ln{L}\bigr)^{-\gamma}\,+\,2(3.01)^{a}L^{-0.25}\,<\,2^{15.6}\bigl(A^{-1}\ln{L}\bigr)^{-\gamma},\quad
\gamma = \frac{1}{2}+\biggl\{\frac{a-1}{2}\biggr\}.
\]
All the above arguments can be applied without any modifications to
the second sum of the theorem. The theorem is proved.

\vspace{0.2cm}

By standard way, we derive the following assertion from the theorems
11 and 12.

\vspace{0.2cm}

\textbf{Theorem 14.} \emph{For a real} $x$ \emph{let quantity}
$\nu_{1} = \nu_{1}(x)$ \emph{denotes the number of ordinates}
$\gamma_{n}$, $N<n\le N+M$ \emph{satisfying the condition}
\[
S(\gamma_{n}+0)\,\le\,\frac{x}{\pi}\sqrt{\frac{L}{2}}.
\]
\emph{Then we have the equation}
\[
\nu_{1}(x)\,=\,M\biggl(\frac{1}{\sqrt{2\pi}}\int_{-\infty}^{\,x}e^{-u^{2}/2}\,+\,\frac{\theta
A}{\sqrt{\ln{L}}}\biggr).
\]
\emph{The quantity} $\nu_{2}(x)$ \emph{of ordinates} $\gamma_{n}$,
$N<n\le N+M$ \emph{with the condition}
\[
S(\gamma_{n}-0)\,\le\,\frac{x}{\pi}\sqrt{\frac{L}{2}}
\]
obeys the similar relation.

Proof of this theorem repeats word\,-for\,-word the proof of theorem
10.

\vspace{0.2cm}

\textbf{Corollary.} \emph{Let} $\Phi(x)$ \emph{be an arbitrary
positive monotonic function that increases unboundedly as} $x\to
+\infty$. \emph{Then the inequalities}
\begin{equation}
\frac{1}{\Phi(n)}\sqrt{\ln\ln
n}\,<\,|S(\gamma_{n}+0)|\,\le\,\Phi(n)\sqrt{\ln\ln n}\label{Lab35}
\end{equation}
\emph{holds for `almost all'} $n$ (\emph{i.e. the quantity of
numbers} $n\le N$ \emph{that do not satisfy inequalities}
(\ref{Lab35}), \emph{is} $o(N)$ \emph{as} $N\to +\infty$).
\emph{This relation holds true if we replace the quantities}
$S(\gamma_{n}+0)$ \emph{in} (\ref{Lab35}) \emph{to}
$S(\gamma_{n}-0)$.

\vspace{0.2cm}

\textbf{Remark.} Let $K\ge K_{0}(\vep)>0$ be a sufficiently large
integer and $N$ runs through all the integers from the interval $(K,
2K]$. Then all the theorems from \S\S 2\,-4 hold true with $M =
[K^{\vep}]$ for each $N$ in this interval, except, may be
$K^{1-0.05\vep}$ values.

\end{document}